\documentclass[12pt]{article}
\usepackage{amsmath,amsthm}
\usepackage{geometry} 
\usepackage[latin1]{inputenc}
\usepackage[english]{babel}
\usepackage{mathrsfs}
\usepackage{array} 
\usepackage{amsmath}
\usepackage{latexsym}
\usepackage{amssymb}

\newtheorem{thm}{Theorem}[section]
\newtheorem{cor}[thm]{Corollary}
\newtheorem{prop}[thm]{Proposition}
\newtheorem{lem}[thm]{Lemma}
\theoremstyle{remark}
\newtheorem{rem}[thm]{Remark}

\theoremstyle{definition}

\newtheorem{ex}[thm]{Example}

\numberwithin{equation}{section}
\numberwithin{thm}{section}

\newcommand{\rmdiv}{{\rm div}}
\newcommand{\psl}{\widetilde{PSL_{2 }(\mathbb{R})}}
\title{Minimal Surfaces in $\psl$}
\author{Rami Younes\\
\small Laboratoire de Mathématiques et Physique Théorique\\[-0.8ex]
\small Université François Rabelais de Tours\\[-0.8ex]
\small \texttt{younes@lmpt.univ-tours.fr}}
\date{{2009}}
\begin{document}
\maketitle
\begin{abstract} We study minimal graphs in the homogeneous Riemannian 3-manifold $\psl$ and we give examples of invariant surfaces. We derive a gradient estimate for solutions of the minimal surface equation in this space and develop the machinery necessary to prove a Jenkins-Serrin type theorem for solutions defined over bounded domains of the hyperbolic plane.
\end{abstract} 
\section{Introduction}
In recent years there has been an increasing interest in the study of minimal and constant mean curvature surfaces in simply connected homogeneous Riemannian 3-manifolds with four dimensional isometry groups. Results in \cite{Uwe}, like the existence of a generalized Hopf-differential or of a Schwarz reflection principle in such manifolds, suggest that these manifolds are the proper setting for studying global properties of minimal and cmc surfaces. The geometries of such manifolds have been classified by Thurston to be either those of the product spaces $\mathbb{S}^2\times\mathbb{R}$ and $\mathbb{H}^2\times\mathbb{R}$, the Heisenberg group $Nil(3)$, or the fiber spaces Berger sphere and $\psl$(see $\cite{scott}$).

Certain aspects of the theory of minimal and cmc surfaces in  $\mathbb{S}^2\times\mathbb{R}$, $\mathbb{H}^2\times\mathbb{R}$ and $Nil(3)$ have been studied for example in, \cite{rosen}, \cite{nelli}, \cite{fmp}, \cite{Dan} and \cite{Uwe} among others.\\In this paper, we study minimal graphs in $\psl$, known to be a Riemannian fibration over the hyperbolic plane, and we obtain a Jenkins-Serrin type theorem for such graphs over convex bounded domains in the hyperbolic plane.We emphasize that $\psl$ is not a product space and so one should ask what is meant by graph in such a space.

A graph in $\psl$ will be the image of a section of the Riemannian submersion $\pi:\psl\to\mathbb{H}^2$.  
A Jenkins-Serrin type theorem gives necessary and sufficient conditions for the solvability of the Dirichlet problem for the minimal surface equation allowing infinite boundary values, prescribed on arcs of the boundary of a convex bounded domain in $\mathbb{H}^2$, and continuous data on the rest of the boundary. However, boundary arcs of a bounded domain in $\mathbb{H}^2$, where a solution of the minimal surface equation in $\psl$ admits infinite values, have to be geodesics (see section 7).

Then more precisely, let $\Omega$ be a convex bounded domain in $\mathbb{H}^2$ whose boundary consists of (open) geodesic arcs $A_1,.., A_n, B_1,..., B_m$, together with their end points and convex open arcs $C_1, C_2,...,C_s$. We suppose that no two geodesics $A_i$ and no two geodesics $B_i$ have a common end point. We give necessary and sufficient conditions for the existence of a minimal section $s:\Omega\to\psl$  of the Riemannian submersion, taking values $+\infty$ on the arcs  $A_1,... A_n$, $-\infty$ on the arcs $B_1,..., B_m$ and arbitrary prescribed continuous data on the arcs $C_1,...,C_s$.\\
For a simple closed geodesic polygon $\mathcal{P}$, whose vertices are chosen from among the endpoints of the segments $A_i$ and the segments $B_i$, let $\alpha$ and $\beta$ be, respectively, the total $\mathbb{H}^2$-length of the geodesics $A_i$ and the total $\mathbb{H}^2$-length of the geodesics $B_i$ which are part of $\mathcal{P}$. Let $\gamma$ be the perimeter of $\mathcal{P}$. Note that in the case $\{C_s\}=\emptyset$, $\mathcal{P}$ could be the whole boundary of $\Omega$.

We have the following
\begin{thm}\label{jsmainthm} If the family of arcs $\{C_s\}$ is non empty, then there exists a unique section of the bundle $\pi:\psl\to\mathbb{H}^2$ defined in $\Omega$ and taking the boundary values $+\infty$ on the geodesics $A_i$, the value $-\infty$ on the geodesics $B_i$ and arbitrary continuous data $f_s$ on $C_s$ if and only if \begin{center}
$2\alpha<\gamma$ and $2\beta<\gamma$
\end{center}for each polygon $\mathcal{P}$ chosen as above.\vspace{.5cm}\\ If the family of arcs $\{C_s\}$ is empty,the condition on the polygons $\mathcal{P}$ is the same except that in the case when $\mathcal{P}$ is the entire boundary of $\Omega$ then the condition is $\alpha=\beta$. Moreover, uniqueness is up to additive constants.
\end{thm}
In $\mathbb{R}^3$ this theorem corresponds to that of Jenkins and Serrin proved in \cite{js}, and in $\mathbb{H}^2\times\mathbb{R}$ a corresponding result was obtained in \cite{nelli} by Nelli and Rosenberg. In their paper, Jenkins and Serrin make use of the a priori estimates for solutions of the minimal surface equation, proved in \cite{serrin}, to obtain a compactness principle for sequences of solutions and to study limit behavior of monotone sequences of solutions. They also make use of the Scherk surface as a barrier, which is fundamental to most of the results. The techniques developed by Serrin in \cite{serrin} were
adapted in \cite{nelli} to show a Jenkins-Serrin type theorem in $\mathbb{H}^2\times\mathbb{R}$. To obtain a priori gradient estimates for solutions of the minimal surface equation and to prove a compactness principle, we adapt a result in Spruck's \cite{sprk} and we construct explicit barriers adequate to our space.

The paper is organized as follows: in section $2$ we give a model for $\psl$, compute its metric in the coordinates and give the expression of its Levi-Civita connection. We then characterize the isometry group of $\psl$ based on ideas from \cite{Bona} and \cite{scott}. We show that this group is generated by the lifts of isometries of $\mathbb{H}^2$ and translations along the fibers.

In section $3$ we derive the minimal surface equation in $\psl$ and furnish examples of minimal graphs invariant under actions of one-parameter groups of isometries generated by lifts of isometries of $\mathbb{H}^2$. The rest of the paper is dedicated to develop the machinery necessary to prove our Jenkins-Serrin type theorem where we follow the main lines in \cite{js}.

In section $4$ we prove an estimate for the gradient of a solution of the minimal surface equation which implies a compactness principle for sequences of solutions of the minimal surface equation in $\psl$ uniformly bounded on compacts of a bounded open subset of $\mathbb{H}^2$.

In section $5$ we prove the existence of a solution of the Dirichlet problem for the minimal surface equation in $\psl$ in a convex bounded open subset of $\mathbb{H}^2$ with boundary data having  possibly a finite number of discontinuities.

 In sections $6$ and $7$ we prove a series of lemmas and propositions which will serve as machinery to prove our Jenkins-Serrin type theorem. Once this machinery is established, the lines of proof are similar to that of the corresponding Jenkins-Serrin theorem in \cite{js} and the reader will be referred to that paper for further details.
\section{The space $\widetilde{PSL_{2 }(\mathbb{R})}$}
The 3-dimensional Lie group of $2\times2$ real matrices of determinant $1$ is denoted $SL_{2}(\mathbb{R})$. The quotient Lie group $SL_{2}(\mathbb{R})$/$\{±I_d\}$ is denoted $PSL_{2 }(\mathbb{R})$ and its universal covering $ \widetilde{PSL_{2 }(\mathbb{R})}$. Of course $ \widetilde{PSL_{2 }(\mathbb{R})}$ is a Lie group itself and so admits left invariant metrics. For our purposes, it will be convenient to introduce a model for $\widetilde{PSL_{2 }(\mathbb{R})}$ and write down explicitly the metric that interests us. In fact we shall show that $ \widetilde{PSL_{2 }(\mathbb{R})}$ is a Riemannian fibration over the hyperbolic plane, the reader can refer to \cite{scott}.
\begin{rem}{A homogeneous simply connected 3-manifold $M$ with a 4-dimensional isometry group, is a Riemannian fibration over a 2-dimensional space form, and whose fibers are geodesics tangent to a unitary Killing field, say $\xi$. These manifolds are classified, up to isometries, by the curvature $\kappa$ of the fibration base and the bundle curvature $\tau$. The number $\tau$ is such that $\overline{\nabla}_X\xi=\tau X\times\xi$, for any vector field $X$ ($\overline{\nabla}$ is the Levi-Civita connection of $M$). As we shall see in what follows, $\psl$ belongs to this class of manifolds  and that the parameters $\kappa$ and $\tau$ have the values $-1$ and $-\frac{1}{2}$ respectively.}\end{rem}
\subsection{A model for $\widetilde{PSL_{2 }(\mathbb{R})}$}
It is known that the group of orientation preserving isometries of the hyperbolic plane $\mathbb{H}^{2}$ is $PSL_2(\mathbb{R})$. Let $U\mathbb{H}^2$ denote the unit tangent bundle of $\mathbb{H}^2$, $i.e.$ the submanifold of $T \mathbb{H}^2$ consisting of tangent vectors of unit length. It is easy to see that $PSL_2(\mathbb{R})$ acts transitively on $U\mathbb{H}^2$ and the stabilizer of each point under this action is trivial. This allows us to identify $PSL_{2 }(\mathbb{R})$ and $U\mathbb{H}^2$  and consequently $\widetilde{PSL_{2 }(\mathbb{R})}$ and $\widetilde{U\mathbb{H}^2}$.
\\The submanifold $U\mathbb{H}^2$ is diffeomorphically a trivial circle bundle over $\mathbb{H}^2$, meaning that $U\mathbb{H}^2 \simeq \mathbb{H}^2\times\mathbb{S}^1$. This implies that $\widetilde{PSL_{2 }(\mathbb{R})}\simeq \mathbb{H}^2\times\mathbb{R}$ again from a diffeomorphic point of view.
\subsection{Metric on $\psl$}
A Riemannian metric on a  manifold $M$ induces a natural metric on the tangent bundle $TM$. We explain how this is generally done and we fix some terminology on the way, the reader can refer to \cite{carmo}. Let $(p,v) \in TM$ and $V$ a tangent vector to $TM$ at $(p,v)$. Choose a curve $\alpha:t\to(p(t),v(t))$ with $p(0)=p, v(0)=v$ and $V=\alpha'(0)$. Define \begin{equation*}
\Vert V\Vert^2_{(p,v)}=\Vert d\pi(V)\Vert^2_{p} + \Vert \frac{Dv}{dt}(0)\Vert^2_{ p},
\end{equation*}
where $\pi: TM\to M$ is the bundle projection and $\displaystyle \frac{D}{dt}$ is the covariant derivative along the curve $t\to p(t)$. The value of $\Vert V\Vert_{(p,v)}$ is independent of the choice of the curve $\alpha$.

A vector at $(p,v)\in TM$ which is orthogonal to the fiber $\pi^ {-1}(p)\simeq T_pM$ is said to be horizontal, and one which is tangent to the fiber is said to be vertical. We identify the vertical tangent space in $T_{(p,v)}(TM)$ to $T_pM$. We have
\begin{enumerate}
	\item[(i)]$\|V\|_{(p,v)}=\|V\|_p$ if $V$ is vertical
	\item[(ii)] $\|V\|_{(p,v)}=\|d\pi(V)\|_p$ if $V$ is horizontal. 
\end{enumerate} 
  Horizontal tangent spaces have the same dimension as tangent spaces to $M$ which implies, together with the identity $(ii)$, that $d\pi$ induces isometries between horizontal tangent spaces and spaces tangent to $M$, $i.e.$, \begin{center}
$ d\pi:TM\to M$
\end{center} is a Riemannian submersion.

Now the metric on $\mathbb{H}^2$ induces a metric on $T\mathbb{H}^2$ which restricts to a metric on $U\mathbb{H}^2$. So we have a metric on $PSL_2(\mathbb{R})$ which lifts to a metric on its universal covering $\widetilde{PSL_{2 }(\mathbb{R})}$. The fact that $PSL_2(\mathbb{R})$ acts on $U\mathbb{H}^2$ by isometries implies that the metric induced on $PSL_2(\mathbb{R})$ is left invariant. This metric lifts obviously to a left invariant metric on $\widetilde{PSL_{2 }(\mathbb{R})}$.

To see that $\psl$ is a Riemannian fibration over $\mathbb{H}^2$ note that the fibres  of $U\mathbb{H}^2$ are 1-dimensional, hence horizontal tangent spaces to $U\mathbb{H}^2$ coincide with those of $T\mathbb{H}^2$ and $\pi$ restricts to a Riemannian submersion on $U\mathbb{H}^2$. As $\widetilde{U\mathbb{H}^2}$ and $U\mathbb{H}^2$ are locally isometric we deduce that $\pi$ induces a Riemannian submersion on $\widetilde{U\mathbb{H}^2}$ onto $\mathbb{H}^2$. The metric on $\widetilde{PSL_{2 }(\mathbb{R})}$ being left invariant (hence complete) we have  $\widetilde{PSL_{2 }(\mathbb{R})}$  a complete homogeneous simply connected Riemannian manifold.

At this point we have given a model for $\widetilde{PSL_{2 }(\mathbb{R})}$ and assigned it a metric. We next express this metric in coordinates, the reader can refer to \cite{Dan}. Let $(x,y)\to\xi(x,y)$ be a conformal parametrization of $\mathbb{H}^2$ and let $\lambda$ be the conformal factor so that the metric of $\mathbb{H}^2$, in these coordinates, is $\lambda^2(dx^2+dy^2)$. As $v\in U\mathbb{H}^2$ is identified with its base point and the angle $\theta$ it makes with $\partial_x$ we have the following local parametrization of $U\mathbb{H}^2$  \begin{center} $\displaystyle (x,y,\theta)\to(\xi(x,y),\frac{1}{\lambda}(\cos\theta \partial_x +\sin\theta \partial_y)).$ \end{center}Let $V$ be a tangent vector to $\widetilde{PSL_{2 }(\mathbb{R})}$ at a point $(p,v)$ and let $\alpha:t\to (p(t),v(t))$ be a curve passing through $(p,v)$ at $t=0$ and tangent to $V$ over there. We write $p(t)=(x(t),y(t))$ and $\displaystyle v(t)=\frac{1}{\lambda}\Big(\cos\theta(t)\partial_x+\sin\theta(t)\partial_y)\Big)$. Using properties of the covariant derivative along the curve $t\to p(t)$ we compute 
\begin{align*}\displaystyle
\frac{Dv}{dt}=& -\frac{\lambda'}{\lambda^2}(\cos\theta\partial_x + \sin\theta\partial_y)+\frac{\theta'}{\lambda}(-\sin\theta\partial_x + \cos\theta\partial_y)\\&+\frac{1}{\lambda}(\cos\theta\nabla_{p'(0)}\partial_x + \sin\theta\nabla_{p'(0)}\partial_y),
\end{align*}
with $\lambda'=x'\lambda_x + y'\lambda_y, p'(0)=x'\partial_x + y'\partial_y, \nabla_{p'(0)}\partial_x=x'\nabla_{\partial_x}\partial_x + y'\nabla_{\partial_y}\partial_x$ and $\nabla_{p'(0)}\partial_y=x'\nabla_{\partial_x}\partial_y+y'\nabla_{\partial_y}\partial_y$. The Christoffel symbols for the metric $\lambda^2(dx^2+dy^2)$ on $\mathbb{H}^2$ are
\begin{eqnarray*} \Gamma_{11}^1=-\Gamma_{22}^1=\Gamma_{12}^2=\Gamma_{21}^2=\frac{\lambda_x}{\lambda}\\  -\Gamma_{11}^2=\Gamma_{22}^2=\Gamma_{12}^1=\Gamma_{21}^1=\frac{\lambda_y}{\lambda}.\end{eqnarray*}
\\We finally obtain
 \begin{equation*}
\frac{Dv}{dt}=\frac{1}{\lambda^2}(\lambda\theta'+y'\lambda_x-x'\lambda_y)(\cos\theta\partial_y - \sin\theta\partial_x).
\end{equation*}
Thus \begin{center} $\|V\|^2_{(p,v)}=\lambda^2(x'^2+y'^2)+\frac{1}{\lambda^2}(\lambda\theta'+y'\lambda_x-x'\lambda_y)^2$.\end{center} Setting $z=\theta$ on the universal covering we get the following expression  for the metric on $\widetilde{PSL_{2 }(\mathbb{R})}$ :\begin{equation*}
ds^2=\lambda^2(dx^2+dy^2)+(-\frac{\lambda_y}{\lambda}dx+\frac{\lambda_x}{\lambda}dy+dz)^2.
\end{equation*}
\begin{rem}{We can see that in our model the fibers are the vertical lines and that a unitary vector field tangent to the fibers is $\xi=\partial_z$. We can also see that translations along the fibers $(x,y,z)\to(x,y,z+a)$ are isometries generated by $\xi$. Thus the fibers are the trajectories of a unit Killing field and so are geodesics.}\end{rem}
\subsection{An orthonormal frame on $\psl$} Let $\{e_1,e_2\}$ be the orthonormal frame on $\mathbb{H}^2$ with $e_1=\lambda^{-1}\partial_x$ and $e_2=\lambda^{-1}\partial_y$ and let $E_3$ be the vector field on $\psl$ whose expression in coordinates is $\xi$. Denote by $E_1$ and $E_2$ the horizontals lifts to $\psl$ of $e_1$ and $e_2$, $i.e.$, \begin{center}
$d\pi(E_i)=e_i$ and $\langle E_i,E_3\rangle=0, 1\leqslant i\leqslant2$.
\end{center}
We remark that $d\pi(\partial x)=\partial x$ and $d\pi(\partial y)=\partial y$, then a simple computation gives the expression of $E_i$ in coordinates, 
\begin{center}
$\displaystyle E_1=\frac{1}{\lambda}\partial_x+\frac{\lambda_y}{\lambda^2}\partial_z, \; E_2=\frac{1}{\lambda}\partial_y-\frac{\lambda_x}{\lambda^2}\partial_z \;\rm{and}\;E_3=\partial_z.$
\end{center}
In what follows let $\tilde{X}$ denote the horizontal lift to $\psl$ of a vector field $X$ on $\mathbb{H}^2$; recall that $\overline{\nabla}_{\tilde{X}}\tilde{Y}=\widetilde{\nabla_XY}+\frac{1}{2}\lbrack\tilde{X},\tilde{Y}\rbrack^v$ for vector fields $X,Y$ on $\mathbb{H}^2$. Then the Riemannian connection of $\psl$ is calculated in the basis $\{E_i\}$ as follows:
 \begin{center}
$\displaystyle\overline{\nabla}_{E_1}E_1=\widetilde{\nabla_{e_1}e_1}=-\frac{\lambda_y}{\lambda^2}E_2, \;\overline{\nabla}_{E_2}E_2=\widetilde{\nabla_{e_2}e_2}=-\frac{\lambda_x}{\lambda^2}E_1.$
\end{center}
 As $E_3$ is a unitary killing field we have, for $\displaystyle1\leqslant i\leqslant3,$ \begin{center}
$\langle \overline{\nabla}_{E_3}E_3,E_i\rangle=-\langle \overline{\nabla}_{E_i}E_3,E_3\rangle=0,$
\end{center} hence, \begin{center}
$\displaystyle \overline{\nabla}_{E_3}E_3=0$.
\end{center} 
For $i,j\in\{1,2\}$ we have, \begin{center}
$\displaystyle\langle \overline{\nabla}_{E_j}E_i,E_j\rangle=-\langle \overline{\nabla}_{E_j}E_j,E_i\rangle$ and $\langle \overline{\nabla}_{E_j}E_i,E_i\rangle=0,$
\end{center}
\begin{center}
$\displaystyle2\langle \overline{\nabla}_{E_i}E_j,E_3\rangle=\langle\lbrack E_i,E_j\rbrack,E_3\rangle-\langle\lbrack E_i,E_3\rbrack,E_3\rangle-\langle\lbrack E_j,E_3\rbrack,E_3 \rangle,$\end{center}
and\begin{center}
$\lbrack E_i,E_3\rbrack=0.$ 
\end{center}
A direct computation of $\lbrack E_1,E_2\rbrack$ gives\begin{center}
$\displaystyle\lbrack E_1,E_2\rbrack=\frac{\lambda_y}{\lambda^2}E_1-\frac{\lambda_x}{\lambda^2}E_2+\Lambda E_3$\end{center}with \begin{center}$\displaystyle\Lambda=\frac{\lambda_x^2+\lambda_y^2}{\lambda^4}-\frac{\lambda_{xx}+\lambda_{yy}}{\lambda^3}=-\frac{\Delta \log\lambda}{\lambda^2}.$\end{center}The last term of the equality is known to be the expression of the curvature, of $\mathbb{H}^2$ in this case, in terms of the conformal factor in isothermal parameters. Therefore, $\Lambda=-1$ and
\begin{center}$\displaystyle\langle \overline{\nabla}_{E_1}E_2,E_3\rangle=-\langle \overline{\nabla}_{E_2}E_1,E_3\rangle=-\frac{1}{2}.$\end{center}
We thus obtain \begin{center}$\displaystyle\overline{\nabla}_{E_1}E_2=\frac{\lambda_y}{\lambda^2}E_1-\frac{1}{2}E_3$\\$\displaystyle\overline{\nabla}_{E_2}E_1=\frac{\lambda_x}{\lambda^2}E_2+\frac{1}{2}E_3.$\end{center}
Moreover the facts that for $1\leqslant i\leqslant2$,
\begin{center}$\displaystyle\lbrack E_i,E_3\rbrack=0,\;\langle \overline{\nabla}_{E_3}E_i,E_i\rangle=0,$\\$\displaystyle\langle \overline{\nabla}_{E_3}E_i,E_3\rangle=-\langle \overline{\nabla}_{E_3}E_3,E_i\rangle=0,$
$\displaystyle\langle \overline{\nabla}_{E_3}E_1,E_2\rangle=\langle \overline{\nabla}_{E_1}E_3,E_2\rangle=-\langle \overline{\nabla}_{E_1}E_2,E_3\rangle=\frac{1}{2}$\\
$\displaystyle\langle \overline{\nabla}_{E_3}E_2,E_1\rangle=\langle \overline{\nabla}_{E_2}E_3,E_1\rangle=-\langle \overline{\nabla}_{E_2}E_1,E_3\rangle=-\frac{1}{2}$\\
\end{center}
conclude that
 \begin{align*}
  \overline{\nabla}_{E_3}E_1&= \overline{\nabla}_{E_1}E_3=\frac{1}{2}E_2,\\
   \overline{\nabla}_{E_3}E_2&=\overline{\nabla}_{E_2}E_3=-\frac{1}{2}E_1.
\end{align*}
We resume our computation \begin{center}$\displaystyle\overline{\nabla}_{E_1}E_1=-\frac{\lambda_y}{\lambda^2}E_2,\;\overline{\nabla}_{E_2}E_2=-\frac{\lambda_x}{\lambda^2}E_1,$\\$\displaystyle\overline{\nabla}_{E_3}E_3=0,$\\$\displaystyle\overline{\nabla}_{E_1}E_2=\frac{\lambda_y}{\lambda^2}E_1-\frac{1}{2}E_3,$\\$\displaystyle\overline{\nabla}_{E_2}E_1=\frac{\lambda_x}{\lambda^2}E_2+\frac{1}{2}E_3,$\\$\displaystyle\overline{\nabla}_{E_3}E_1= \overline{\nabla}_{E_1}E_3=\frac{1}{2}E_2,$\\$\displaystyle \overline{\nabla}_{E_3}E_2=\overline{\nabla}_{E_2}E_3=-\frac{1}{2}E_1.$\end{center}
\begin{rem}{The equation $\overline{\nabla}_{E_3}E_3=0$ is the geodesic equation for vertical fibers.}\end{rem}
\begin{rem}{The fact that $[E_1,E_2]$ is not horizontal implies that the horizontal plane field generated by $E_1$ and $E_2$ is not integrable, meaning that there exists no horizontal surfaces in $\psl$.}\end{rem}
\subsection{Isometries of $\widetilde{PSL_{2 }(\mathbb{R})}$}
It is known that $\psl$ has a four dimensional isometry group. See \cite{scott} for example, a standard reference on the geometries of 3-manifolds. However in \cite{scott} this group is  characterized using Lie group theory. In what follows is what the author of this paper found a worth while simplified geometric characterization of this group based on ideas from \cite{scott} and $\cite{Bona}$.

The metric induced on the tangent bundle $TM$ of a Riemannian manifold M is intrinsic enough that it is respected by the lifts of isometries of M to $TM$. In fact, each map $f\in C^\infty (M,M)$ lifts to a map $df\in C^\infty(TM,TM)$ such that\\$\displaystyle df(p,v)=(f(p),d_pf(v))$. When $f$ is an isometry, $df$ induces isometries on tangent spaces of $TM$. This can be easily seen as follows. Let $(p,v)\in TM$ and $V\in T_{(p,v)} (TM)$ and choose a curve $\alpha(t)=(p(t),v(t))$ in $TM$ such that $\alpha(0)=(p,v)$ and $\alpha'(0)=V$. We have,
\begin{equation*}
\|d_v(df)V\|_{(f(p),d_pf(v)}^2=\|d_pf(p'(0))\|_{f(p)}^2+\|\frac{Ddf(v)}{dt}(0)\|_{f(p)}^2,
\end{equation*}
where $\displaystyle \frac{D}{dt}$ is the covariant derivative along the curve $\beta(t)=df(\alpha(t))$. As $d_pf$ is an isometry and \begin{center}$\displaystyle \frac{Ddf(v)}{dt}=df(\frac{Dv}{dt})$\end{center}it follows directly that
\begin{equation*}
\|d_{(p,v)}(df)V\|_{(f(p),d_pf(v))}=\|V\|_{(p,v)},
\end{equation*}  proving our claim.\\
In particular, the isometry group of $\widetilde{PSL_{2 }(\mathbb{R})}$ contains the lifts of the isometries of $\mathbb{H}^2$. We note also that vertical translations along the fibers are isometries of $\widetilde{PSL_{2 }(\mathbb{R})}$. These isometries read in coordinates as $(x,y,z)\to(x,y,z+a)$.
So the isometry group of $\widetilde{PSL_{2 }(\mathbb{R})}$ contains the group $G$ generated by the lifts of isometries of $\mathbb{H}^2$ and vertical translations. In fact, we shall show that $G$ contains all the isometries of $\widetilde{PSL_{2 }(\mathbb{R})}$. We begin with  proving the following proposition found in $\cite{Bona}$.
\begin{prop}
The sectional curvature along a plane $P\subset T_{(p,v)}(\widetilde{PSL_{2 }(\mathbb{R})})$ is maximal  when P contains the line $L_{(p,v)}$, the line tangent to the fiber at $(p,v)$, and is minimal when P is orthogonal to $L_{(p,v)}$.
\end{prop}
\begin{proof}Let P be a plane generated by two orthonormal vectors X and Y. Then the sectional curvature along $P$ is given by $\left\langle \overline{R}(X,Y)X,Y \right\rangle$, where $\overline{R}$ is the curvature tensor of $\widetilde{PSL_{2 }(\mathbb{R})}$. We have $\displaystyle \langle R(X,Y)X,Y \rangle=\frac{-7}{4}+2(\left\langle X,\xi \right\rangle^2+\left\langle Y,\xi \right\rangle^2)$ (see \cite{Dan}, proposition 2.1). As $\xi$ is unitary we have $\left\langle X,\xi \right\rangle^2+\left\langle Y,\xi \right\rangle^2\leqslant1$. So the sectional curvature will be maximal when $\left\langle X,\xi \right\rangle^2+\left\langle Y,\xi \right\rangle^2=1$, and this is possible only when $\left\langle \xi,Z \right\rangle=0$ for any vector $Z$ such that $\{X,Y,Z\}$ forms an orthonormal basis of the tangent space to $\psl$ at $(p,v)$. This means that the sectional curvature will be maximal when $\xi \in P$, $i.e.$ when $P$ contains the vertical line tangent to the fiber. Similarly we show that the sectional curvature is minimal when $P$ is orthogonal to the vertical line tangent to the fiber.\end{proof}
We next show that isometries of $\widetilde{PSL_{2 }(\mathbb{R})}$ are fiber preserving. The proposition above implies that the differential of an isometry $\varphi$ sends $L_{(p,v)}$ to $L_{\varphi(p,v)}$. This follows from the fact that the differential of an isometry will send two planes  along which the sectional curvature is maximal, to two planes along which the curvature is maximal. As the fiber $\pi^{-1}(p)$, tangent at the point $(p,v)$ to $L_{(p,v)}$, is a geodesic, its image under $\varphi$ is the geodesic tangent to the line $L_{\varphi(p,v)}$ at the point $\varphi(p,v)$. The fiber through $\varphi(p,v)$ is a geodesic tangent to the former line at $\varphi(p,v)$, so we conclude that it is the geodesic in question. We have then the following
\begin{prop}
The isometries of $\widetilde{PSL_{2 }(\mathbb{R})}$ are fiber preserving, $i.e.$ the images by an isometry of two points lying on the same fiber belong to the same fiber.
\end{prop}
This property will allow each isometry of $\widetilde{PSL_{2 }(\mathbb{R})}$ to induce an isometry on $\mathbb{H}^2$ the following manner,
\begin{lem}
Every isometry $\varphi$ on $\widetilde{PSL_{2 }(\mathbb{R})}$ induces an isometry $f$ on $\mathbb{H}^2$ such that $f\circ \pi=\pi\circ \varphi$.\end{lem}
\begin{proof} The equation $f\circ \pi=\pi\ \circ \varphi$ defines $f$ the obvious way as $\varphi$ is fiber preserving. For a vector $v\in T_p\mathbb{H}^2$ such that $v= d_{(p,v)}\pi(V)$, $V$ is the horizontal lift of $v$, we have $d_pf(v)=d_pf(d_{(p,v)}\pi(V))=d_{\varphi(p,v)}\pi(d_{(p,v)}\varphi (V))$. As $V$ is horizontal and $\varphi$ is an isometry we have $d_{(p,v)}\pi(V)$ also horizontal. The fact that $\pi$  is a Riemannian submersion concludes that $f$ is indeed an isometry.\end{proof}
We proceed to show the following technical lemma found in \cite{scott}, which will aid giving the finishing touch to our characterization of isometries of $\widetilde{PSL_{2 }(\mathbb{R})}$.
\begin{lem}
Fix a point $(p,v) \in \widetilde{PSL_{2 }(\mathbb{R})}$. We may compose any isometry $\alpha$ of $\widetilde{PSL_{2 }(\mathbb{R})}$ with isometries lying in $G$ to obtain an isometry $\beta$ which fixes $(p,v)$ and whose differential at $(p,v)$ is the identity on the horizontal tangent plane at $(p,v)$.\end{lem}
\begin{proof}Let $f$ be the isometry induced by $\alpha$ on $\mathbb{H}^2$. We compose $\alpha$ with a vertical translation sending $\alpha(p,v)$ to $(f(p),d_pf(v))$ to obtain an isometry $\alpha'$ of $\widetilde{PSL_{2 }(\mathbb{R})}$. Let $df^{-1}$ denote the lift of $f^{-1}$ to $\widetilde{PSL_{2 }(\mathbb{R})}$ and set $\beta=df^{-1} \circ \alpha'$. This is an isometry of $\widetilde{PSL_{2 }(\mathbb{R})}$ fixing $(p,v)$ and leaving each horizontal vector at $(p,v)$ invariant. In fact, for a horizontal vector $V$ at $(p,v)$ we have \begin{center}$d_{(p,v)}\beta(V)=d_{(f(p),d_pf(v))}df^{-1}(d_{(p,v)}\alpha(V)).$ \end{center}
We denote the restriction of $d\pi$ to horizontal tangent planes by $d\pi_{\circ}$ and we set $w=d_p\pi(V)$, so we have\begin{center} $d_{(p,v)}\alpha(V)=d_pf(w)$ \rm{and} $d_{(f(p),d_pf(v))}df^{-1}(d_{(p,v)}\alpha(V))=d_p\pi_{\circ}^{-1}(w)=V$.\end{center}We used the fact that $d_{(p,v)}dg(V)=d_p\pi_{\circ}^{-1}(d_pg(w))$, for any lift $dg$ of an isometry $g$ of $\mathbb{H}^2$.\end{proof}
At this point it is easy to prove our claim that $G$ contains all the isometries of $\widetilde{PSL_{2 }(\mathbb{R})}$. Let $\varphi$ be an isometry and $(p,v)$ a point of $\widetilde{PSL_{2 }(\mathbb{R})}$. We compose $\varphi$ with isometries in $G$ and we obtain an isometry $\psi$ which fixes $(p,v)$, and whose differential at $(p,v)$ is the identity on the horizontal plane at $(p,v)$. Consequently $\psi$ leaves invariant the fiber through $(p,v)$ as it is fiber preserving.

Let $\ell$ be a piecewise geodesic loop in $\mathbb{H}^2$ based at $p$ with non-trivial holonomy and $\tilde{\ell}$ be its horizontal lift to $\widetilde{PSL_{2 }(\mathbb{R})}$ starting at $(p,v)$. Let $(p,w)$ denote the other end of $\tilde{\ell}$.
Now $ \psi(\tilde{\ell})$ is piecewise geodesic since so is $\tilde{\ell}$ and as $\psi$ is an isometry (see Remark $6$ below). Since $\psi$ fixes $(p,v)$ and the horizontal plane over there we deduce that $\psi(\tilde{\ell})$ passes through $(p,v)$ and has the same horizontal tangent vector as $\tilde{\ell}$ there. Hence $\psi(\tilde{\ell})$ equals $\tilde{\ell}$ and in particular $\psi$ must fix $(p,w)$.\\
As $\psi$ is an isometry and the points $(p,v)$ and $(p,w)$ are distinct, due to non trivial holonomy of the geodesic loop based at $p$ below in $\mathbb{H}^2$, it follows that $\psi$ fixes each point of the fiber through $(p,v)$. Then $\psi$ is an isometry which fixes a point and whose differential over there is the identity.This implies that $\psi$ leaves invariant geodesics through $(p,v)$. As our manifold is complete we can join $(p,v)$ to any other point by a geodesic. Being an isometry $\psi$ fixes each point of these geodesics and so $\psi$ is the identity. This allows us to deduce that $\varphi$ is a composition of elements of $G$.

We resume the result in the following,
\begin{thm}\label{psliso}
The isometry group of $\widetilde{PSL_{2 }(\mathbb{R})}$ is generated by the lifts of the isometries of $\mathbb{H}^2$ together with the vertical translations along the fibers.
\end{thm}
\begin{rem}{Theorem \ref{psliso} implies that the isometry group of $\psl$ is four dimensional and contains no orientation reversing isometries.}\end{rem} 
\begin{rem}{
Assume that $\gamma:t\to\gamma(t)$ is a geodesic in $\mathbb{H}^2$ starting at a point $p$. We can lift $\gamma$ to a horizontal geodesic in $\psl$, one whose velocity vector at each point is horizontal, starting at any point $(p,v)$ in the fiber above $p$.
Fix such a point $(p,v)$ and let $v(t)$ be the parallel transport of $v$ along $\gamma$.
The curve $\bar{\gamma}:t\to(\gamma(t),v(t))$ starts at $(p,v)$. The fact that $v(t)$ is parallel implies that $\bar{\gamma}$ is horizontal. To show that $\bar{\gamma}$ a geodesic we suppose to the contrary that it is not. We choose convex neighborhoods $W\subset \psl$ of  $(p,v)$  and $U\subset \mathbb{H}^2$ of $p$ such that $\pi(W)=U$. Take two points $Q_1=(q_1,w_1)$ and $Q_2=(q_2,w_2)$ in $\bar{\gamma}\cap W$, joined by an arc $\bar{\alpha}$ such that $L(\bar{\alpha})< L(\bar{\gamma})=L(\gamma)$. In $\mathbb{H}^2$, $\gamma$ is a minimizing geodesic joining $q_1$ and $q_2$. The arc $\alpha=\pi(\bar{\alpha})$ verifies $L(\alpha)\leqslant L(\bar{\alpha})$, which contradicts the fact that $\gamma$ is length minimizing \rm(see \cite{carmo}, p.79). }\end{rem}
\section{Minimal graphs in $\psl$}
We fix our model of $\psl$ as $\mathbb{H}^2\times\mathbb{R}$ endowed with the metric \begin{equation*}
ds^2=\lambda^2(dx^2+dy^2)+(-\frac{\lambda_y}{\lambda}dx+\frac{\lambda_x}{\lambda}dy+dz)^2,
\end{equation*} as described above.\vspace{.5cm}\\We denote by $S_\circ \subset \psl$ the surface defined by $z=0$. We identify a domain $\Omega\subset \mathbb{H}^2$ and its lift to $S_\circ$. We define the graph $\Sigma(u)$ of $u\in C^0(\bar{\Omega})$ on $\Omega$ as \begin{equation*}
\Sigma(u)=\{(x,y,u(x,y))\in \psl | (x,y)\in \Omega\}.
\end{equation*}
These graphs are basically images of sections of the bundle projection \begin{center}
$\pi:\psl\to \mathbb{H}^2,$
\end{center} 
$i.e.$ images of maps $s:\Omega \subset \mathbb{H}^2\to \psl$ with $\pi\circ s= I_{\mathbb{H}^2}$.
For such a map let $u(x,y)$ be the signed distance from the lift of $(x,y)\in\mathbb{H}^2$, the point of $\psl$ whose coordinates are $(x,y,0)$, to $s(x,y)\in \pi^{-1}(x,y)$ along the geodesic fiber through $(x,y,0)$. The fibers here being oriented positively by $\xi$. This function $u$ defined by $s$ defines a graph, in the sense of the above definition, which is the image of $s$. Clearly, each function $u\in C^0(\bar{\Omega}), \Omega\subset\mathbb{H}^2$, defines a section of the bundle projection.\vspace{.5cm}\\For a smooth function $u$ set $F(x,y,z)=z-u(x,y)$ so that $\Sigma(u)=F^{-1}(0)$. As F is smooth we will have \begin{center}
$\displaystyle\eta=\frac{\overline{\nabla}F}{|\overline{\nabla}F|}$
\end{center} a unit normal field to $\Sigma(u)$.\\A simple computation shows that
 \begin{center}
$\displaystyle \overline{\nabla}F=(\frac{\lambda_y}{\lambda^2}-\frac{u_x}{\lambda})E_1+(-\frac{\lambda_x}{\lambda^2}-\frac{u_y}{\lambda})E_2+E_3.$
\end{center}
Set\begin{center} $\displaystyle\alpha=\frac{\lambda_y}{\lambda^2}-\frac{u_x}{\lambda}, \beta=-\frac{\lambda_x}{\lambda^2}-\frac{u_y}{\lambda}$ and $W=| \overline{\nabla}F|=\sqrt{1+\alpha^2+\beta^2}$,\end{center}so that \begin{center} $\displaystyle\eta=\frac{\alpha}{W}E_1+\frac{\beta}{W}E_2+\frac{1}{W}E_3.$\vspace{.5cm}\end{center}
We parameterize the graph of a smooth function $u$ by \begin{center}$(x,y)\to\phi(x,y)=(x,y,u(x,y)),$ \end{center}with $(x,y)\in\Omega$ the domain of definition of $u$. It is easy to see that for the metric on $\psl$ we have
\begin{center}
$\displaystyle\langle\phi_x,\phi_x\rangle=\lambda^2(1+\alpha^2)$, $\displaystyle\langle\phi_x,\phi_y\rangle=\lambda^2\alpha\beta$, $\displaystyle\langle\phi_y,\phi_y\rangle=\lambda^2(1+\beta^2),$
\end{center}giving the metric induced on the graph \begin{center}
$\displaystyle g=\lambda^2\big((1+\alpha^2)dx^2+\alpha\beta dxdy+\alpha\beta dydx+(1+\beta^2)dy^2\big).$
\end{center}To calculate the mean curvature $H$ of $\Sigma(u)$, with respect to the upwards pointing normal $\eta$, choose $v_1,v_2\in T(\psl)$ so that $\big\{v_1,v_2,\eta\big\}$ is an orthonormal basis of $T(\psl)$. As $\eta$ is a unitary field we have $\langle\overline{\nabla}_\eta\eta,\eta\rangle=0$ and
\begin{align*} \displaystyle
2H=&-\sum_{1}^{2} \langle \overline{\nabla}_{v_i}\eta,v_i\rangle\\
      =&-\sum_{1}^{2} \langle \overline{\nabla}_{v_i}\eta,v_i\rangle-\langle \overline{\nabla}_\eta\eta,\eta\rangle\\
      =&-\rmdiv{(\eta)}.
\end{align*}
Therefore $\displaystyle 2H=-\rmdiv\left(\frac{\overline{\nabla}\textsl{F}}{| \overline{\nabla}\textsl{F}|}\right),$ where div and $\overline{\nabla}$ denote respectively the divergence and the Levi-Civita connection in $\psl$.\\
Since $E_1$ and $E_2$ are the horizontal lifts of $e_1$ and $e_2$ , the facts that $\overline{\nabla}_{E_3}E_3=0$ and that $\pi$ is a Riemannian submersion allow us to write 
\begin{align*}\displaystyle
\rmdiv\left(\frac{\alpha}{W}E_1+\frac{\beta}{W}E_2\right)
    =&\sum_{1}^{2}\Big\langle\overline{\nabla}_{E_i} \left(\frac{\alpha}{W}E_1+\frac{\beta}{W}E_2\right),E_i\Big\rangle_{\psl}\\
    =&\sum_{1}^{2}\Big\langle\nabla_{e_i}d\pi\left(\frac{\alpha}{W}E_1+\frac{\beta}{W}E_2\right),e_i\Big\rangle_{\mathbb{H}^2}\\
    =&\rmdiv _{\mathbb{H}^2}\Big(\frac{\alpha}{\lambda W}\partial_x+\frac{\beta}{\lambda W}\partial y\Big).
\end{align*}
Since $E_3$ is a Killing field we have $\rmdiv(E_3)=0$, and \begin{center}
$\displaystyle \rmdiv \Big(\frac{1}{W}E_3\Big)=\Big\langle \overline{\nabla}\Big(\frac{1}{W}\Big),E_3\Big\rangle+\frac{\rmdiv(E_3)}{W}=\frac{\partial }{\partial z}\Big(\frac{1}{W}\Big)=0$.
\end{center}
Therefore \begin{center}
$\displaystyle2H= \rmdiv_{\mathbb{H}^2}\left(\frac{\alpha}{\lambda W}\partial_x+\frac{\beta}{\lambda W}\partial_y\right)=\rmdiv_{\mathbb{H}^2}\Big(d\pi(\eta)\Big).$
\end{center}
We also have,\begin{center}
$\displaystyle2H=\frac{1}{\lambda^2}\rmdiv_{\mathbb{R}^2}\left(\frac{\lambda\alpha}{W}\partial_x+\frac{\lambda\beta}{W}\partial_y\right)$,
\end{center} as $\displaystyle \rmdiv_{\mathbb{H}^2}(X)=\frac{1}{\lambda^2}\rmdiv_{\mathbb{R}^2}(\lambda^2X)$ for any vector field $X$ on $\mathbb{H}^2$.
The equation of a minimal graph is then \begin{equation}\label{mse}
\rmdiv_{\mathbb{R}^2}\left(\frac{\lambda\alpha}{W}\partial_x+\frac{\lambda\beta}{W}\partial_y\right)=0.
\end{equation}
\subsection{Examples of minimal surfaces and minimal graphs in $\psl$}
In this section we find minimal graphs invariant under the action of the one parameter groups of isometries of $\psl$ generated by the lifts of rotations, parabolic and hyperbolic isometries of $\mathbb{H}^2$. We also determine the minimal surfaces invariant under translation along the fibers. 
\begin{ex} \label{liftsofgeo}\rm{
Let $\gamma$ be a geodesic of $\mathbb{H}^2$. The vertical cylinder over $\gamma$, $\mathcal{C}_\gamma=\pi^{-1}(\gamma) \subset \psl$, is a minimal surface and this can be seen as follows:
Let $T$ and $\eta$ be respectively a unit tangent field and a unit normal field to $\gamma$, and let $\widetilde{T}$ and $\widetilde{\eta}$ be their corresponding horizontal lifts to $\psl$. We then have $\{\widetilde{T},E_3\}$ an orthonormal basis on $\mathcal{C}_\gamma$ and $\widetilde{\eta}$ a unit normal to $\mathcal{C}_\gamma$. The mean curvature of $\mathcal{C}_\gamma$ at a point $v$ is then 
\begin{align*}
2H&=-\Big\langle\bar{\nabla}_{\widetilde{T}}\widetilde{\eta},\widetilde{T}\Big\rangle-\Big\langle \bar{\nabla}_{E_3}\widetilde{\eta},E_3\Big\rangle-\Big\langle \bar{\nabla}_{\widetilde{\eta}}\widetilde{\eta},\widetilde{\eta}\Big\rangle\\
&=\Big\langle \widetilde{\nabla^{\mathbb{H}^2}_TT},\widetilde{\eta}\Big\rangle=\Big\langle \nabla^{\mathbb{H}^2}_TT,\eta\Big\rangle\\
&=\rm{the\:geodesic\;curvature\;of\;\gamma\;at\;the\;point\;\pi(v)},
\end{align*}
and as $\gamma$ is a geodesic we deduce that $H=0$, and the cylinder $\mathcal{C}_\gamma$ is thus minimal. We notice that these minimal surfaces are invariant under vertical translations and they are in fact the only ones. A minimal surface invariant under vertical translations is $\pi^{-1}(\gamma)$, where $\gamma$ is a curve of $\mathbb{H}^2$. The geodesic curvature of $\gamma$ is shown again by the above computation to be zero and hence $\gamma$ is geodesic.
}\end{ex}
\begin{ex}\label{thebarrier}
\rm{ The 1-parameter group of isometries of $\mathbb{H}^2$, given in the half plane model of $\mathbb{H}^2$ by $(x,y)\to(\epsilon x,\epsilon y)$, induces a 1-parameter group of isometries on $\psl$. In our model of $\psl$ these isometries read as $(x,y,z)\to(\epsilon x, \epsilon y, z)$. A minimal graph invariant by this group of isometries is that of a solution $u$ of $(\ref{mse})$
verifying $u(r,\theta)=u(\theta)$, $(r,\theta)$ are polar coordinates on the upper half plane . Here we have \begin{center}
$\displaystyle \lambda=\frac{1}{y}$,\;$\alpha=-yu_x-1$ and $\beta=-yu_y,\;y>0$.
\end{center}
Let $\displaystyle \omega:=W^2=1+\alpha^2+\beta^2$, equation $(\ref{mse})$ then implies \begin{equation}
\label{E1}
\omega\Big( \frac{\partial}{\partial x}(\lambda \alpha)+\frac{\partial}{\partial y}(\lambda\beta) \Big)-\frac{\lambda}{2}(\alpha\omega_x+\beta\omega_y)=0.
\end{equation}
An invariant solution $u$ verifies
\begin{center}
$\displaystyle u_x =\frac{\partial\theta}{\partial x}u_{\theta} =-\frac{\sin\theta}{r}u_\theta,$\\
$\displaystyle u_y =\frac{\partial\theta}{\partial y} u_{\theta}=\frac{\cos\theta}{r}u_\theta,$\\
$\displaystyle \omega=2-2\sin^2\theta u_{\theta}+\sin^2\theta u_{\theta}^2,$\\
$\displaystyle u_{xx}=\Big(\frac{\partial\theta}{\partial x}\Big)^2 u_{\theta\theta}+\frac{\partial^2\theta}{\partial^2x}u_{\theta}=\frac{\sin^2\theta}{r^2}u_{\theta\theta}+2\frac{\sin\theta \cos\theta}{r^2}u_\theta,$\\
$\displaystyle u_{yy}=\Big(\frac{\partial\theta}{\partial y}\Big)^2u_{\theta\theta}+\frac{\partial^2\theta}{\partial^2y}u_{\theta}=\frac{\cos^2\theta}{r^2}u_{\theta\theta}-2\frac{\sin\theta\cos\theta}{r^2}u_{\theta}.$
\end{center}
Equation $(\ref{E1})$ implies that 
\begin{equation}\label{E3}\displaystyle
\omega u_{\theta \theta}-\frac{1}{2}\omega_{\theta}(u_{\theta}-1)=0
\end{equation}
from which we deduce that either
\begin{enumerate}
	\item[(i)]  $u_\theta=1$, or
	\item[(ii)] $\displaystyle 2\frac{u_{\theta\theta}}{u_\theta-1}=\frac{\omega_\theta}{\omega}$
\end{enumerate}
which is equivalent to \begin{equation*}\displaystyle\frac{(u_\theta-1)^2}{\omega}=C,\;C\geq0.\end{equation*}
The cases $(i)$ and $(ii)$ are resumed in
\begin{equation*}\displaystyle
(1-C\sin^2\theta )(u_{\theta}^2-2u_{\theta})=2C-1,C\geq0.
\end{equation*}
For $0\leqslant C<1$, this first integral defines a 1-parameter family of graphs over the hyperbolic plane, given up to an additive constant by
\begin{center}
$\displaystyle u(r,\theta)=u(\theta)=\pm \sqrt{C}\int_{0}^{\theta} \frac{\sqrt{1+\cos^2\theta}}{\sqrt{1-C\sin^2\theta}}d\theta +\theta,\;0<\theta<\pi.$\end{center}For example, when $C=0$ we obtain up to vertical translations, half a (euclidean) Helicoid over the hyperbolic plane.\\When $\displaystyle C=\frac{1}{2}$ the above solutions simplify to $u(r,\theta)=\theta\pm\theta+constant$. So on the one hand we obtain up to vertical translations, half a Helicoid stretched in the vertical direction. It is the surface over the hyperbolic plane obtained by rotating, in euclidean terms, the $x-$axis about the $z-$axis, and translating it vertically twice as fast. On the other hand we obtain translates of the plane $\{z=0\}$ as invariant minimal surfaces which correspond to the solutions $u(r,\theta)=constant$.\\
For $C=1$, we obtain solutions defined in the first and the second quadrants of the hyperbolic plane. The solutions are\begin{center}$\displaystyle u(r,\theta)=u(\theta)= \int_{0}^{\theta} \frac{\sqrt{1+\cos^2\theta}}{\cos\theta}d\theta +\theta$\\$\displaystyle\Big(= -\int_{0}^{\theta} \frac{\sqrt{1+\cos^2\theta}}{\cos\theta}d\theta +\theta\;\rm{respectively}\Big),\;0<\theta<\frac{\pi}{2},$\end{center}defined in the first quadrant and taking values $0$ on the positive $x-$axis and $+\infty$ ($-\infty$ respectively) on the $y-$axis. On the other hand the solutions\begin{center}$\displaystyle u(r,\theta)=u(\theta)= \int_{\frac{\pi}{2}}^{\theta} \frac{\sqrt{1+\cos^2\theta}}{\cos\theta}d\theta +\theta$\\$\displaystyle\Big(= -\int_{\frac{\pi}{2}}^{\theta} \frac{\sqrt{1+\cos^2\theta}}{\cos\theta}d\theta +\theta,\;\rm{respectively}\Big),\frac{\pi}{2}<\theta<\pi,$\end{center} defined in the second quadrant and taking values $+\infty$ ($-\infty$ respectively) on the $y-$axis and $0$ on the negative $x-$axis. The solutions obtained so far define complete minimal graphs.\\
For $C>1$, the equation $1-C\sin^2\theta=0$ has two solutions, say $\theta_1$ and $\theta_2=\pi-\theta_1$, in $]0,\pi[$ such that $\theta_1<\frac{\pi}{2}<\theta_2$. The first integral defines a one-parameter family of disconnected graphs defined in the region $\{0<\theta<\theta_1\}\bigcup\{\theta_2<\theta<\pi\}$. We have, up to additive constants, the solutions\begin{center}$\displaystyle u(r,\theta)=u(\theta)=\sqrt{C}\int_{0}^{\theta} \frac{\sqrt{1+\cos^2\theta}}{\sqrt{1-C\sin^2\theta}}d\theta +\theta$\\$\displaystyle\Big(=-\sqrt{C}\int_{0}^{\theta}\frac{\sqrt{1+\cos^2\theta}}{\sqrt{1-C\sin^2\theta}}d\theta +\theta\;\rm{respectively}\Big),\;0<\theta<\theta_1,$\end{center}and\begin{center}$\displaystyle u(r,\theta)=u(\theta)=\sqrt{C}\int_{\theta}^{\pi} \frac{\sqrt{1+\cos^2\theta}}{\sqrt{1-C\sin^2\theta}}d\theta +\theta$\\$\displaystyle\Big(=-\sqrt{C}\int_{\theta}^{\pi}\frac{\sqrt{1+\cos^2\theta}}{\sqrt{1-C\sin^2\theta}}d\theta +\theta\;\rm{respectively}\Big),\;\theta_2<\theta<\pi.$\end{center}One can see easily that the solutions have finite values over the lines $\theta=\theta_1$ and $\theta=\theta_2$ and admit vertical tangent planes over there. However, the solutions obtained for these values of the parameter $C$ do not define complete minimal graphs.\\ One obtains complete minimal surfaces above the region $\{0<\theta\leq\theta_1\}$ for example, when one considers unions of graphs $u(r,\theta)$ above that region. We consider the graphs obtained for both factors $\pm\sqrt{C}$ of the integral in the above expression of $u$ and translate them vertically to take values $\theta_1$ over $\theta=\theta_1$.To see that the union defines a regular surface above $\theta=\theta_1$ we simply show that $\theta$ is a smooth function of $z$ near $z=\theta_1$.\\We have $z=u(\theta)$ which implies that the derivatives of $\theta$ with respect to $z$ are given by\begin{center}$\displaystyle\frac{\partial u}{\partial\theta}=\frac{1}{\theta'}$ \;and\; $\displaystyle\frac{\partial^2 u}{\partial \theta^2}=-\frac{\theta''}{\theta'^3}.$\end{center} We compute $\omega$ and $\omega_{\theta}$ in terms of $\theta$ and its derivatives then substitute in $(\ref{E3})$ to obtain after necessary simplifications,\begin{equation} \label{invhom}\displaystyle\theta''(2-\sin^2\theta)+\sin\theta\cos\theta(\theta'-1)(2\theta'-1)=0.\end{equation}As the graphs $u(r,\theta)$ admit vertical tangent planes at the points $z=\theta_1$, $\theta$ defines a $C^1$-function of $z$ and the equation $(\ref{invhom})$ shows then that $\theta$ is in fact smooth.$\hfill\Box$
}\end{ex}
\begin{ex} \rm{
Consider the disc model for the hyperbolic plane. Rotations, in euclidean terms, about the center of the disc are isometries of $\mathbb{H}^2$. The lifts of these isometries to $\psl$, seen in our model, are euclidean screw motions. The image of a point $(x,y,z)$ is obtained by rotating the $(x,y)$ part around the $z-$axis then translating it along the $z-$axis by the same amount. We can then compose the lift of a rotation on $\mathbb{H}^2$ with a translation along a vertical fiber to obtain an isometry of $\psl$ which is rotation about the fiber. So we have a 1-parameter group of isometries of $\psl$ which are, in our model, rotations about the $z-$axis.\\A minimal graph invariant by this group is that of a solution $u$ of $(\ref{E1})$ verifying $u(r,\theta)=u(r),\;(r,\theta)$ polar coordinates on the disc.\\Here we have,
 \begin{center}
$\displaystyle \lambda=\frac{1}{1-\frac{x^2+y^2}{4}}$,\;$\displaystyle \alpha=-\frac{u_x}{\lambda}+\frac{y}{2}$ and $\displaystyle \beta=-\frac{u_y}{\lambda}-\frac{x}{2}, \; x^2+y^2<4$.
\end{center}
An invariant solution verifies
\begin{center}
$\displaystyle u_x =\frac{\partial r}{\partial x}u_{r} =\frac{x}{r}u_r,$\\
$\displaystyle u_y =\frac{\partial r}{\partial y} u_{r}=\frac{y}{r}u_r,$\\
$\displaystyle \omega=1+\frac{r^2}{4}+\frac{1}{\lambda^2}u_r^2,$\\
$\displaystyle u_{xx}=\Big(\frac{\partial r}{\partial x}\Big)^2 u_{rr}+\frac{\partial^2r}{\partial^2x}u_{r}=\frac{x^2}{r^2}u_{r r}+\frac{y^2}{r^3}u_r,$\\
$\displaystyle u_{yy}=\Big(\frac{\partial r}{\partial y}\Big)^2u_{r r}+\frac{\partial^2r}{\partial^2y}u_{r}=\frac{y^2}{r^2}u_{r r}+\frac{x^2}{r^2}u_{r}.$
\end{center}
Equation $(\ref{E1})$ implies that 
\begin{center}
$\displaystyle \omega(u_{r r}+\frac{1}{r}u_r)-\frac{1}{2}u_r\omega_r=0,$
\end{center}
from which we deduce that either
\begin{enumerate}
	\item[(i)]$u\equiv constant$, or
	\item[(ii)]$\displaystyle 2\frac{u_{rr}}{u_r}+\frac{2}{r}=\frac{\omega_r}{\omega}$
\end{enumerate}
which is equivalent to\begin{center}$\displaystyle r^2u_r^2=C\omega,\;C>0.$\end{center}This implies that \begin{center}$\displaystyle u_r=\pm2\sqrt{\frac{r^2+4}{Cr^2-(r^2-4)^2}}$,\end{center} with $C>0$ and $0<r_\circ<r<2$,
$\displaystyle r_\circ=\sqrt{\frac{8+C-\sqrt{(8+C)^2-64}}{2}}.$\vspace{.5cm}\\Remark that $u_r(r_{\circ})=\pm\infty$ and that the solutions are either increasing or decreasing in $r$. In a fashion similar to that in the above example, we show that the union of the graphs corresponding to both values of $u_r$ and taking the value $0$ at $r_\circ$ define a regular surface. Then this first integral defines up to vertical translations, a family of minimal surfaces of catenary type. As $C$ varies in $]0,+\infty[$ the asymptotic angles at infinity between the members of the family and the cylinder $\partial\mathbb{H}^2\times\mathbb{R}$ assume all the values in $]0,\pi[$.\\Remark also that up to a vertical translation, when $C\rightarrow+\infty$, $r_\circ\rightarrow0$ and the limit surface is the doubly covered hyperbolic plane(identified with $z=0$). When $C\rightarrow0$, $r_{\circ}\rightarrow2$ and up to a vertical translation the family degenerates to the circle at infinity $\partial\mathbb{H}^2$ (doubly covered).
}\end{ex}
\begin{ex}\label{paratype}
\rm{ The 1-parameter group of isometries of $\mathbb{H}^2$, given in the half plane model of $\mathbb{H}^2$ by $(x,y)\to( x+a, y)$, induces a 1-parameter group of isometries on $\psl$. In our model of $\psl$ these isometries read as $(x,y,z)\to(x+a, y,z) $.\\A minimal graph invariant by this group of isometries is that of a solution $u$ of $(\ref{mse})$
verifying $u(x,y)=u(y),\;y>0$.\\We have  
\begin{center}
$\displaystyle \lambda=\frac{1}{y},\;\alpha=-1,\;\beta=-yu_y$ and $\omega=2+y^2u_y^2,$
\end{center}
and so equation $(\ref{E1})$ implies that
\begin{center}
$\displaystyle \omega u_{yy}-\frac{1}{2}u_y\omega_y=0.$
\end{center}
We deduce that either $u\equiv constant$, or $\displaystyle u_y=\pm\frac{\sqrt{2}}{\sqrt{C^2-y^2}},\;C>0.$ This equation defines up to additive constants, surfaces symmetric (in euclidean terms) with respect to $\{z=0\}$. These surfaces are the union of the two graphs \begin{center}
$\displaystyle u(x,y)=\pm\sqrt{2}\arcsin\Big(\frac{y}{C}\Big)\mp\frac{\sqrt{2}\pi}{2}$\end{center}over the region $\{0<y\leq C\}$. As $C\rightarrow+\infty$ the limit surface is $\{z=0\}$ (doubly covered).
}\end{ex}$\hfill\Box$
\begin{rem}{
There exists no compact complete minimal surface in $\psl$. For otherwise, if such a surface $\Sigma$ exists, we may then translate down any minimal surface $z=constant$ not intersecting $\Sigma$ until there is a first contact point. This implies that the two surfaces are tangent and one above the other. By the maximum principle $\Sigma$ will be equal to a surface $z=constant$. This is a contradiction as the surfaces $z=constant$ are not compact.
}\end{rem}
\section{Gradient Estimates}
We will next prove an estimate for the gradient of a solution $u:\Omega\subset\mathbb{H}^2\to\mathbb{R}$ of $(\ref{mse})$ following the lines of proof of Theorem 1.1 in \cite{sprk}, which will be fundamental for proving later results. For this aim we will need the following formulae which hold for surfaces in 3-manifolds and in particular for surfaces $\Sigma\subset\psl$:
\begin{eqnarray}
\label{e5}|\nabla_\Sigma f|^2&=&|\overline{\nabla} \tilde{f}|^2-\langle\overline{\nabla} \tilde{f},\eta\rangle^2\\
\label{e6}\Delta_\Sigma f&=&2\langle\overline{\nabla} \tilde{f},\eta\rangle H+\Delta\tilde{f}-\langle\overline{\nabla}_\eta\overline{\nabla}\tilde{f},\eta\rangle\\
\label{e7}\Delta_\Sigma g(f)&=&g'(f)\Delta_\Sigma f+g''(f)|\nabla_\Sigma f|^2,\end{eqnarray}
where $f$ is a function defined on $\Sigma$, or the restriction to $\Sigma$ of a function \begin{center}$\tilde{f}:\psl\rightarrow\mathbb{R}^2$,\end{center}$H$ is the mean curvature of $\Sigma$ and $\eta$ a unit normal field on $\Sigma$. We will also need the following fact, if $X:M\rightarrow N$ is a constant mean curvature isometric immersion of a surface $M$ in a 3-manifold $N$, and if $\eta$ is a unit normal field to $M$ and $\xi$ a Killing field on $N$ then the function $n=\langle \eta,\xi\rangle$ verifies the following equation\begin{equation}\label{e8}
\Delta_\Sigma n=-(|A|^2+Ric(\eta))n,
\end{equation} where $|A|$ is the norm of the second fundamental form of $M$ and $Ric$ is the Ricci curvature of $N$.\\For minimal graphs in $\psl$, $\displaystyle n=\frac{1}{W}$, so that\begin{equation}\label{stability}\displaystyle \Delta_\Sigma \frac{1}{W}=-(|A|^2+Ric(\eta))\frac{1}{W}\end{equation}with\begin{equation*}
Ric(\eta)=-\frac{3}{2}+\frac{2}{W^2},
\end{equation*}which we compute using the equations of proposition 2.1 in \cite{Dan}.
\begin{rem}{\label{stablegraphs}Equation $(\ref{stability})$ implies that minimal graphs in $\psl$ are stable. This follows directly from the definition of stability of a minimal surface and Theorem 1 in \cite{fc}.}\end{rem}\rm{}

We finally note that a function $\phi:\Omega\to\mathbb{R}$ lifts as a section of $\pi$ to a function on $\psl$, whose restriction to $\Sigma$ will be also denoted by $\phi$. Then using $(\ref{e5})$ we obtain \begin{center}
$\displaystyle|\nabla_\Sigma\phi|_\Sigma^2=\frac{1}{\lambda^2W^2}\big((\phi_x^2+\phi_y^2)+(\beta\phi_x-\alpha\phi_y)^2\big)$\end{center}
which implies that\begin{equation}\label{e10}
\displaystyle|\nabla_\Sigma\phi|_\Sigma^2\geqslant\frac{1}{W^2}|D\phi|_{\mathbb{H}^2}^2.\end{equation} 
\begin{thm}\label{gradientestimates}
Let $u$ be a non-negative solution of the minimal surface equation $(\ref{mse})$ in a bounded domain $\Omega\subset \mathbb{H}^2$. Then at each point $p\in\Omega$ we have\begin{center}
$\displaystyle W(p)\leqslant C$\end{center}where $C$ is a positive constant which depends only on $u(p)$, the distance of $p$ to $\partial\Omega$ and on bounds of $\lambda$ and its derivatives on $\Omega$.
\end{thm}
\begin{proof} We fix a point $p\in \Omega$. We introduce the function $f=\mu(x)W$ on a geodesic ball $B_\rho(p)\subset\Omega\subset\mathbb{H}^2$, for which we will derive a maximum principle by computing $\Delta_\Sigma f$. The function $\mu$ is to be defined. We have \begin{align*}
  \Delta_\Sigma f&=W\Delta_\Sigma\mu+2\langle\nabla_\Sigma W,\nabla_\Sigma\mu\rangle+\mu\Delta_\Sigma W \\
                    &=W\Delta_\Sigma\mu+\frac{2}{W}\big(\langle\nabla_\Sigma W,\nabla_\Sigma f\rangle-\mu|\nabla_\Sigma W|^2\big)+\mu\Delta_\Sigma W.
\end{align*}
We then obtain \begin{center}
$\displaystyle\Delta_\Sigma f-\frac{2}{W}\langle\nabla_\Sigma W,\nabla_\Sigma f\rangle=\mu\big(\Delta_ \Sigma W-\frac{2}{W}|\nabla_\Sigma W|^2\big)+W\Delta_\Sigma\mu.$
\end{center}
However from $(\ref{e7})$ we get 
\begin{center}
$\displaystyle\Delta_\Sigma \frac{1}{W}=-\frac{1}{W^2}\Delta_\Sigma W+\frac{2}{W^3}|\nabla_\Sigma W|^2,$
\end{center}and $(\ref{e8})$ then implies that
\begin{center}
$\displaystyle\Delta_\Sigma W-\frac{2}{W}|\nabla_\Sigma W|^2=(|A|^2+Ric(\eta))W,$
\end{center}
so that \begin{center}
$\displaystyle\Delta_\Sigma W-\frac{2}{W}|\nabla_\Sigma W|^2\geqslant Ric(\eta)W\geqslant-\frac{3}{2} W$
\end{center}
We get 
\begin{center}
$\displaystyle\Delta_\Sigma f-\frac{2}{W}\langle\nabla_\Sigma W,\nabla_\Sigma f\rangle\geqslant W\big(\Delta_\Sigma\mu-\frac{3}{2}\mu\big).$
\end{center}
The idea is to define $\mu$ so that $\displaystyle\Delta_\Sigma\mu-\frac{3}{2}\mu>0.$ \\We set\begin{center}
$\displaystyle\mu(x)=e^{K\phi}-1,$ and $\displaystyle\phi(x)=-\frac{u(x)}{2u_\circ}+1-\Big(\frac{d(x)}{\rho}\Big)^2$
\end{center}
on the ball $B(p,\rho)$, where $u_\circ=u(p)$, $d$ is the geodesic distance from $p$ and $K>0$ a constant to be determined. We next bound $\Delta_\Sigma\mu-\frac{3}{2}\mu$ from below. Using $(\ref{e7})$ we obtain\begin{center}
$\displaystyle\Delta_\Sigma\mu=K e^{K\phi}\Delta_\Sigma\phi+K^2e^{K\phi}|\nabla_\Sigma\phi|^2.$
\end{center}
As $\displaystyle u=h_{|_\Sigma}$, $h=z$ in the given model of $\psl$ and $\Sigma$ minimal, $(\ref{e6})$ implies that\begin{center}
$\Delta_\Sigma\ u=\Delta h-\langle\overline{\nabla}_\eta\overline{\nabla}h,\eta\rangle,$
\end{center}
showing that we can bound $\Delta_\Sigma u$ by a constant independent of $u$. Similarly we bound $\Delta_\Sigma d^2$ which shows that\begin{center}
$\displaystyle\Delta_\Sigma\phi\geqslant-C_1\big(\frac{1}{u_\circ}+\frac{1}{\rho^2}\big), $
\end{center}where $C_1$ is a constant.\vspace{.5cm}
The inequality $(\ref{e10})$ implies that in $\textit{B}_\rho(p)$\begin{center}
$\displaystyle|\nabla_\Sigma\phi|_\Sigma^2\geqslant\frac{1}{W^2}|D\phi|_{\mathbb{H}^2}^2\geqslant\frac{1}{W^2}\Big(\frac{|Du|_{\mathbb{H}^2}^2}{4u^2_\circ}-\frac{2}{u_\circ\rho}|Du|_{\mathbb{H}^2}\Big),$\end{center}which implies that when\begin{center}$\displaystyle|Du|_{\mathbb{H}^2}\geqslant\frac{16u_\circ}{\rho}$\end{center}we have\begin{center}
$\displaystyle|\nabla_\Sigma\phi|_\Sigma^2\geqslant\frac{|Du|_{\mathbb{H}^2}^2}{8u_\circ^2W^2}.$
\end{center}Now as \begin{equation}\label{ine11}
\displaystyle W^2\leqslant1+2|Du|_{\mathbb{H}^2}^2+2\Big(\big(\frac{\lambda_x}{\lambda^2}\big)^2+\big(\frac{\lambda_y}{\lambda^2}\big)^2\Big)
\end{equation} we obtain \begin{center}
$\displaystyle\frac{|Du|_{\mathbb{H}^2}^2}{W^2}\geqslant C_2\frac{|Du|_{\mathbb{H}^2}^2}{1+|Du|_{\mathbb{H}^2}^2},$
\end{center}where $C_2$ is a positive constant which depends only on bounds of $\lambda$ and its derivatives over $\Omega$.
Hence on the set where $\displaystyle|Du|_{\mathbb{H}^2}> max(1,\frac{16u_\circ}{\rho})$ we find \begin{center}
$\displaystyle\Delta_\Sigma\mu-\frac{3}{2}\mu\geqslant C' e^{K\phi}\Big(\frac{C}{u^2_\circ}K^2-\big(\frac{1}{u_\circ}+\frac{1}{\rho^2}\big)K-1\Big),$
\end{center}where $C$ and $C'$ are positive constants which depend on bounds of $\lambda$ and its derivatives.
We next choose \begin{center}
$\displaystyle K>\frac{u^2_\circ}{2C}\Big(\frac{1}{u_\circ}+\frac{1}{\rho^2}+\sqrt{\big(\frac{1}{u_\circ}+\frac{1}{\rho^2}\big)^2+\frac{4C}{u^2_\circ}}\Big)$
\end{center} so that $\displaystyle\Delta_\Sigma\mu-\frac{3}{2}\mu>0$ on the set  $\displaystyle|Du|_{\mathbb{H}^2}> max(1,\frac{16u_\circ}{\rho})$.\\If \begin{center}$|Du|_{\mathbb{H}^2}\leq max(1,\frac{16u_\circ}{\rho})$\end{center} then inequality $(\ref{ine11})$ proves our claim on $W(p)$. Otherwise, we consider the open set \begin{center}$\displaystyle U=\{x\in\textit{B}_\rho(p)/\phi>0, |Du|_{\mathbb{H}^2}>max(1,\frac{16u_\circ}{\rho})\},$\end{center}and note that $p\in U$.Then by the maximum principle, the point $p_\circ$ where $f$ achieves its maximum on $U$ belongs to $\partial U$ with $f(p_\circ)>0$. As $\phi<0$ on $\partial\textit{B}_\rho(p)$ we have $$\partial U\cap\partial\textit{B}_\rho(p)=\emptyset$$ and therefore $$p_\circ\in\{\displaystyle|Du|_{\mathbb{H}^2}=max(1,\frac{16u_\circ}{\rho})\}\cap\{\phi>0\}.$$ Therefore\begin{center}
$\displaystyle f(p)=\mu(p)W(p)\leqslant \mathcal{C}\mu(p_\circ)\sqrt{1+max^2(1,\frac{16u_\circ}{\rho})}$
\end{center}
and $$\displaystyle W(p)\leqslant \mathcal{C}e^{\frac{K}{2}}\sqrt{1+max^2(1,\frac{16u_\circ}{\rho})},$$ where $\mathcal{C}$ is a positive constant which depends only on bounds of $\lambda$ and its derivatives. The proof is completed.\end{proof}
\begin{cor}Let $u$ be a bounded solution of the minimal surface equation $(\ref{mse})$ in a domain $\Omega\subset\mathbb{H}^2$. Then at any point $p\in\Omega$ we have\begin{center}$W(p)\leq\mathcal{C}$\end{center}where $\mathcal{C}$ is a positive constant which depends only on $\displaystyle\max_{\partial\Omega}\left|u\right|$, the distance of $p$ to $\partial\Omega$ and on bounds of $\lambda$ and its derivatives on $\Omega$.
\end{cor}
\begin{thm}\label{uniformbounds}
Let $u$ be a solution of the minimal surface equation $(\ref{mse})$ in $\Omega$ with $W\leqslant M$ at a point $p\in\Omega$. Then there exists $R$, which depends only on $M, u(p)$ and $d(p,\partial\Omega)$, such that $W\leqslant2M$ on D(p,R).
\end{thm}
Proof. We shall derive an estimate on $\|\nabla W\|$, the norm of the $\mathbb{R}^2$-gradient of $W$, from which the bound on $W$ follows readily. The graph of $u$ is parametrized by \begin{center}
$(x,y)\longrightarrow\psi(x,y)=(x,y,u(x,y)),$
\end{center} and a unit normal field to the graph is\begin{center}
$\displaystyle\eta=\frac{\alpha}{W}E_1+\frac{\beta}{W}E_2+\frac{1}{W}E_3.$
\end{center} 
The partial derivatives of $\psi,$\begin{center}
$\psi_x=\partial_x+u_x\partial_z=\lambda E_1-\lambda\alpha E_3$\end{center}
and\begin{center}
$\psi_y=\partial_y+u_x\partial_z=\lambda E_2-\lambda\beta E_3,$
\end{center}
are such that\begin{center}
$\|\psi_x\|_{\psl}^2\leqslant\lambda W$ and $\|\psi_y\|_{\psl}^2\leqslant\lambda W$.
\end{center}
We shall estimate the partial derivatives of $\alpha$ and $\beta$ by applying the Schoen curvature estimate. For this purpose we need to calculate $\|\overline{\nabla}_{\psi_x}\eta\|$,
\begin{align*}
\displaystyle
 \overline{\nabla}_{\psi_x}\eta=&\frac{\partial}{\partial x}\Big(\frac{\alpha}{W}\Big)E_1+\frac{\alpha}{W}(\lambda\overline{\nabla}_{E_1}E_1-\lambda\alpha\overline{\nabla}_{E_3}E_1) \\
    &+\frac{\partial}{\partial x}\Big(\frac{\beta}{W}\Big)E_2+\frac{\beta}{W}(\lambda\overline{\nabla}_{E_1}E_2-\lambda\alpha\overline{\nabla}_{E_3}E_2)\\
    &+\frac{\partial}{\partial x}\Big(\frac{1}{W}\Big)E_3+\frac{1}{W}(\lambda\overline{\nabla}_{E_1}E_3-\lambda\alpha\overline{\nabla}_{E_3}E_3),\\
\end{align*}
so that \begin{center}$\displaystyle \overline{\nabla}_{\psi_x}\eta=U+V$\end{center}with \begin{align*}\displaystyle
U=&\Big(\frac{\partial}{\partial x}\Big(\frac{\alpha}{W}\Big)+\frac{\lambda_y}{\lambda}\frac{\beta}{W}\Big)E_1\\
&+\Big(\frac{\partial}{\partial x}\Big(\frac{\beta}{W}\Big)-\frac{\lambda_y}{\lambda}\frac{\alpha}{W}\Big)E_2\\
&+\frac{\partial}{\partial x}\Big(\frac{1}{W}\Big)E_3.
\end{align*}
and \begin{center}$\displaystyle V=\frac{\lambda\alpha\beta}{2W}E_1+\frac{\lambda(1-\alpha^2)}{2W}E_2-\frac{\lambda\beta}{2W}E_3.$\end{center}It is easy to see that \begin{center}$\displaystyle\|U\|_{\psl}^2\leqslant 2\Big(\|\overline{\nabla}_{\psi_x}\eta\|^2_{\psl}+\|V\|_{\psl}^2\Big).$ \end{center}
We wish to estimate $\|U\|_{\psl}$ as it is the term which contains derivatives of $\alpha$ and $\beta$. We have\begin{align*}
 \frac{\partial}{\partial x}\Big(\frac{\alpha}{W}\Big)&=\frac{1}{W^3}\Big((1+\beta^2)\alpha_x-\alpha\beta\beta_x\Big)\\
 \frac{\partial}{\partial x}\Big(\frac{\beta}{W}\Big)&=\frac{1}{W^3}\Big((1+\alpha^2)\beta_x-\alpha\beta\alpha_x\Big)\\
\frac{\partial}{\partial x}\Big(\frac{1}{W}\Big)&=-\frac{\alpha\alpha_x+\beta\beta_x}{W^3}.
\end{align*}
Therefore,\begin{center}
$\displaystyle \|U\|_{\psl}^2=\frac{1}{W^4}(\alpha_x^2+\beta_x^2)+\Big(\frac{1}{W^2}(\alpha_x\beta-\alpha\beta_x)+\frac{\lambda_y}{\lambda}\Big)^2-\Big(\frac{\lambda_y}{\lambda}\Big)^2\frac{1}{W^2}.$
\end{center}Its easy to see that \begin{center}$\|V\|_{\psl}\leqslant\lambda W$.\end{center}
The shape operator of the graph, which is stable (c.f. Remark $\ref{stablegraphs}$), is $\displaystyle \widetilde{A}\psi_x=-\overline{\nabla}_{\psi_x}\eta$. Schoen's curvature estimate implies that $|\widetilde{A}|\leqslant C$ in a disc about each point on the graph, where $C$ is a constant which depends only on the $\psl$ distance of the point from the boundary of the graph. The inequality\begin{center}
 $\|\widetilde{A}\psi_x\|_{\psl}\leqslant|\widetilde{A}|\|\psi_x\|_{\psl},$
\end{center} implies that at each point $p\in\Omega$
\begin{center}
$\displaystyle \alpha_x^2+\beta_x^2\leqslant \lambda^2CW^6+\lambda^2W^6+\Big(\frac{\lambda_y}{\lambda}\Big)^2W^2,$
\end{center}
and yet \begin{center}
$\displaystyle \alpha_x^2+\beta_x^2\leqslant CW^6,$
\end{center}
$C$ is a constant which depends only on $u(p)$, the distance of $p$ from $\partial\Omega$ and on bounds of $\lambda$ and its derivatives over compacts of $\Omega$.\\Similarly we obtain \begin{center}$\overline{\nabla}_{\psi_y}\eta=U'+V'$\end{center}with
\begin{align*}
\displaystyle
 U'=&\Big(\frac{\partial}{\partial y}\Big(\frac{\alpha}{W}\Big)-\frac{\lambda_x}{\lambda}\frac{\beta}{W}\Big)E_1\\
&+\Big(\frac{\partial}{\partial y}\Big(\frac{\beta}{W}\Big)+\frac{\lambda _x}{\lambda}\frac{\alpha}{W}\Big)E_2\\
&+\frac{\partial}{\partial y}\Big(\frac{1}{W}\Big)E_3,
\end{align*}and\begin{center}$\displaystyle V'=\frac{\lambda(\beta^2-1)}{2W}E_1-\frac{\lambda\alpha\beta}{2W}E_2+\frac{\lambda\alpha}{2W}E_3$\end{center}
The facts \begin{center}$\displaystyle \|U'\|_{\psl}^2=\frac{1}{W^4}(\alpha_y^2+\beta_y^2)+\Big(\frac{1}{W^2}(\alpha\beta_y-\beta\alpha_y)+\frac{\lambda_x}{\lambda}\Big)^2-\Big(\frac{\lambda_x}{\lambda}\Big)^2\frac{1}{W^2},$
\end{center}and \begin{center}$ \| V'\|\leq \lambda W$\end{center}imply that
\begin{center}
$\displaystyle \alpha_y^2+\beta_y^2\leqslant CW^6,$\end{center}$C$ is a constant which depends only on $\Omega$, $u(p)$ and the distance of $p$ from the boundary of $\Omega$.\\Note that $\displaystyle\nabla W=\frac{1}{W}(\alpha\alpha_x+\beta\beta_x,\alpha\alpha_y+\beta\beta_y)$, hence the estimates obtained on the partial derivatives of $\alpha$ and $\beta$ imply that at each point $p\in\Omega,$\begin{center}
$\|\nabla W\|\leqslant CW^3.$
\end{center}This estimate will allow us to conclude our proof. Let $R=\frac{1}{2}d_{\mathbb{R}^2}(p,\partial\Omega)$ and introduce the function $f(r)=W(r,\theta)$ in $D(p,R)\subset\Omega$, where $r$ and $\theta$ are the polar coordinates with origin $p$. We fix $\theta\neq0$ and we remark that $f(0)=W(p)\leq M$ and\begin{center}
$\displaystyle f'(r)=\frac{\partial W}{\partial r}\leqslant\|\nabla W\|\leqslant Cf(r)^3.$
\end{center} 
Integrating this inequality we obtain that $f(r)\leqslant2M$ for $r\in[0,\frac{3}{8M^2C}[,$ which reads into $W$ is bounded by $2M$ on $D(p,min(R,\frac{3}{8M^2C}))$.\vspace{.5cm}$\hfill\Box$\\
The above estimates imply that the first and second derivatives of a solution $u$ at a point $p$, admit bounds which depend only on the value of $u$ at $p$, the distance of $p$ from the boundary and on $\Omega$. Then the classical Ascoli theorem implies the following\\
\textbf{Compactness principle.} Let $(u_n)$ be a uniformly bounded sequence of solutions of the minimal surface equation $(\ref{mse})$ in a domain $\Omega$. Then there exists a subsequence which converges to a solution in $\Omega$, the convergence being uniform on every compact subset of $\Omega$.
\section{Preliminary Existence Theorems }
In what follows $C$ will denote a  rectifiable Jordan curve in $\psl$. Let $\mathcal{D}$ denote the solution of the Plateau problem for $C$ (exists as $\psl$ is homogeneous, see \cite{morrey}), a compact minimal disc with least area, having  $C$ as boundary. It is known that $\mathcal{D}$ has a tangent plane at each interior point, see $\cite{law}$. Let $h$ denote the function defined on $\psl$ whose expression in the model described above is $h=z$ and set $\displaystyle m_C=\min_C(h)$ and $\displaystyle M_C=\max_C(h)$. We suppose that $m<M$ for our curve $C$ for otherwise $\mathcal{D}$ will be a piece of a surface defined by $h=constant$.\\For a curve $\gamma \subset \mathbb{H}^2$, we denote  $\mathcal{C}(\gamma)$ the convex hull of $\gamma$, $i.e.$ the smallest (geodesically) convex subset of $\mathbb{H}^2$ containing $\gamma$ and $R_C=\pi^{-1}(\mathcal{C}\big(\pi(C)\big))$, the region in $\psl$ above the convex hull of the projection of $C$. Note that $C$ is contained in $R_C$. The following proposition corresponds in $\mathbb{R}^3$ to the result that a minimal surface is contained in the convex hull of its boundary.
\begin{prop}\label{convexhull}
The minimal disc $\mathcal{D}$ is contained in $R_C\bigcap \{m_C\leqslant h\leqslant M_C\}$.
\end{prop}
\begin{proof}There exists a minimal disc $\Delta$ defined by $h=constant$ not intersecting $\mathcal{D}$. If $\mathcal{D}$ had an interior point $p$ above (respectively below) all other points of $C$, we would translate $\Delta$ downwards (respectively upwards) along vertical fibers so that $\Delta$ is eventually tangent to $\mathcal{D}$. This is impossible by the maximum principle as we assume $h$ non-constant on $C$. Therefore $\mathcal{D}$ is contained in $\{m_C\leqslant h\leqslant M_C\}$. Similarly we show that $\mathcal{D}$ is contained in $R_C$, except that instead of considering minimal discs $h=constant$, we consider cylinders above geodesics of $\mathbb{H}^2$ and instead of vertical translation we use the fact that these cylinders foliate $\psl$. Note that the interior of $\mathcal{D}$ is strictly contained in the interior of $R_C\bigcap \{m_C\leqslant h\leqslant M_C\}$.\end{proof}
The next proposition asserts the existence of a solution for the Dirichlet's problem for the minimal surface equation in $\psl$, over a convex bounded domain of $\mathbb{H}^2$ with prescribed continuous boundary data.
\begin{prop}[Rado's Lemma in $\psl$]\label{rado} 
If $C$ admits a one-to-one projection onto a convex curve in $\mathbb{H}^2$, then the interior of $\mathcal{D}$ can be obtained as the image of a minimal section of $\pi$.\end{prop}
\begin{proof} Let $C$ be a curve in $\mathbb{H}^2\times\mathbb{R}$, as described above, which has a one-to-one projection onto a convex curve of $\mathbb{H}^2$. We want to prove that the interior of $\mathcal{D}$ is a graph over $\Omega$, the open convex subset of $\mathbb{H}^2$ bounded by the $\pi(C)$. \\
Consider a vertical translate $\mathcal{D'}$ of $\mathcal{D}$, above $\mathcal{D}$, such that $\mathcal{D}\cap \mathcal{D'}=\emptyset$. We suppose that $\mathcal{D}^\circ$ is not a graph, so that there are two distinct points $P$ and $Q$ of $\mathcal{D}^\circ$, say $P$ above $Q$, lying on the same fiber. Let $P'$ and $Q'$ be the corresponding translates of $P$ and $Q$ on $\mathcal{D'}$. We can translate $\mathcal{D'}$ down as to have $P\equiv Q'$. So at one point, when translating $\mathcal{D}'$ down, a translate $\mathcal{D'}$ will have a first point of contact with $\mathcal{D}$ without having $\mathcal{D}\equiv \mathcal{D'}$. By the maximum principle, this point of contact is not interior to both discs. So either the interior of one disc will touch the boundary of the other, or the boundaries of both discs touch at first. However, the above proposition shows that the interior of each disc lies in $R^\circ_C$, and the boundaries lie on $\partial R_C$ as they have convex projections to $\mathbb{H}^2$. So we are left with the only possibility that the first point of contact is a boundary point for both, which is a contradiction for the boundary is projected one-to-one into $S_{\circ}$.\end{proof}
In the next proposition we show that it is possible to claim existence of solutions when boundary data has a finite set of discontinuities. We will first prove the existence of a particular minimal graph which will be of use as a barrier later on.
\begin{lem}\label{particularbarrier}
Let $\mathcal{T}$ be an isosceles geodesic triangle in $\mathbb{H}^2$ with (open) sides $\mathcal{S}_i$ such that $ length(\mathcal{S}_1)=length(\mathcal{S}_2)$, and $c \in \mathbb{R}^*$. Let $\Delta$ denote the open bounded region of $\mathbb{H}^2$ bounded by $\mathcal{T}$. There exists a solution $u$ of the minimal surface equation $(\ref{mse})$ defined in $\Delta\cup\{\mathcal{T}-$ vertices of $\mathcal{S}_i\}$ such that $u=0$ on $\mathcal{S}_1$ and $\mathcal{S}_2$, and $u=c$ on $\mathcal{S}_3$.
\end{lem}
\begin{proof} Consider such a triangle in $S_\circ$ and let $C\subset \psl$ be the Jordan curve formed by $\mathcal{S}_1$, $\mathcal{S}_2$, the translate of $\mathcal{S}_3$ to height $h=c$, and the two fiber segments joining the vertices of $\mathcal{S}_3$ to those of its translate. Let $\Sigma$ be the interior of the solution of the Plateau problem for $C$. We shall show that $\Sigma$ is a graph, thus showing the existence of our minimal section with the desired values on $\partial\mathcal{T}$.\\Assume to the contrary that $\Sigma$ is not a graph, so that there exist two points $P$ and $Q$ of $\Sigma$ lying on the same fiber, say $P$ above $Q$, with $d(P,Q)=d>0$. Let  $f_\epsilon$ be a family of isometries of $\mathbb{H}^2$ converging to the identity in $\rm{C^1}-$topology, such that  \begin{center}
$f_\epsilon(\mathcal{S}_i)\bigcap \mathcal{C}(\mathcal{T})=\emptyset,\;i=1,\,2.$ 
\end{center}
Let $\tilde{f_\epsilon}$ denote the lift of $\displaystyle f_\epsilon$ to $\psl$ as explained in $2.4$, and $\Sigma_{\epsilon,t}=\tilde{f_\epsilon}(\Sigma)+(0,0,t),\;c.t>0$. For $\left|t\right|\geq d$ and $\epsilon$ small enough, we have\begin{center}
$\displaystyle \Sigma_{\epsilon,t}\cap C=\emptyset,\;\rm{and}\;\partial\Sigma_{\epsilon,t}\cap\Sigma=\emptyset.$
\end{center}
We suppose, without loss of generality, that $c>0$ and we remark that for 
$\epsilon$ small enough we'll have \begin{center}
$\displaystyle \|\tilde{f_\epsilon}-Id_{\psl}\|_\infty<\frac{d}{2}.$
\end{center}
To see the former equality we remark that the boundary of $\displaystyle \Sigma_{\epsilon,d}$ is composed of the arcs $C_{\epsilon,i}=\tilde{f_\epsilon}(\mathcal{S}_i)+(0,0,d)$, plus the fiber segments joining the extremities of
$C_{\epsilon,1}$ to $C_{\epsilon,3}$ and $C_{\epsilon,2}$ to $C_{\epsilon,3}$.
We can see that \begin{center}
$h_{|C_{\epsilon,i}}>\frac{d}{2}$, $(i=1,\:2)$, and $h_{|C_{\epsilon,3}}>c+\frac{d}{2}.$
\end{center}
Then these inequalities show that for $\epsilon$ small enough $\displaystyle \Sigma_{\epsilon,d}$ is above $z=\frac{d}{2}$ and hence \begin{center}
$\displaystyle \Sigma_{\epsilon,d}\cap \mathcal{S}_i=\emptyset.$
\end{center}Moreover, $\displaystyle \Sigma_{\epsilon,d}$ lies in $\pi^{-1}(f_\epsilon(\mathcal{T}))$ by proposition $\ref{convexhull}$, so that \begin{center}
$\displaystyle \Sigma_{\epsilon,d}\cap \mathcal{C}_{\mathcal{S}_3}=\emptyset,$ \end{center}where $\mathcal{C}_{\mathcal{S}_3}$ is the cylinder above $\mathcal{S}_3$, completing the proof that $\displaystyle \Sigma_{\epsilon,d}\cap C=\emptyset$.\\To show that $\partial\Sigma_{\epsilon,d}\cap\Sigma=\emptyset$, we first need to remark that $\Sigma \subset\pi^{-1}(\mathcal{T})$. This implies that $\Sigma$ cannot intersect but possibly $C_{\epsilon,3}$ of $\partial\Sigma_{\epsilon,d}$. However the fact that $z_{|C_{\epsilon,3}}>c+\frac{d}{2}$ shows no intersection in this case either as $\Sigma$ is below $z=c$.\\Therefore, \begin{center}
$\partial\Sigma_{\epsilon,d}\cap\Sigma=\emptyset$.
\end{center}
Now the maximum principle implies that for $\epsilon$ small enough we have \begin{center}
$\Sigma_{\epsilon,d}\cap\Sigma=\emptyset.$
\end{center}If we let $\epsilon\to 0$ we shall thus obtain that the limit surface, $\Sigma_d=\Sigma+(0,0,d)$, tangent to $\Sigma$ at $P\in \Sigma$. By the maximum principle the two surfaces should be equal; a contradiction. Therefore, $\Sigma$ is a graph as was claimed.\end{proof}
We now extend the result of proposition $\ref{rado}$ to include Jordan curves containing finitely many vertical fiber segments.
\begin{prop}\label{radogen}
Let $\Omega$ be a bounded convex domain in $\mathbb{H}^2$ and consider a finite set of boundary points of $\Omega$. Let $C$ denote the remaining boundary of $\Omega$, which consists of a finite number of open arcs. Then there exists a solution of the minimal surface equation in $\Omega$ taking preassigned bounded continuous data on the arcs $C$.\end{prop}
\begin{proof} Let $f$ be the bounded continuous data on $C$ and $f_n$ a bounded sequence of continuous functions on $\partial\Omega$ which converges uniformly to $f$ on compacts of $C$. Let $u_n$ be the solution of the minimal surface equation in $\Omega$ with boundary values $f_n$. Proposition $\ref{convexhull}$ implies that the sequence $u_n$ is uniformly bounded on compact sets of $\Omega$, and hence by the compactness principle admits a subsequence which converges to a solution $u$ in $\Omega$.\\
The function $u$ takes the values $f$ on $C$ as shown below using a standard barrier technique. Indeed, there exist barriers at each point of $C$, i.e., at each point $P$ of $C$ and for each pair of positive numbers $K$ and $\delta$, there exist a neighborhood $V$ of $P$ and a non-negative solution $v$ in $V\cap\Omega$ such that
\begin{enumerate}
	\item[(i)]$V\cap\Omega$ is contained in the geodesic disc of radius $\delta$ about $P$,
	\item[(ii)] $v\geqslant K$ on $\partial V\cap \Omega$,
	\item[(iii)]$v=0$ at $P$.
\end{enumerate}
We may take $V$ to be an isosceles triangle, having its equal sides intersecting in $\Omega$ and tangent to $\partial\Omega$ at $P$ on its third side, and $v$ the solution in this triangle which takes values $K$ on the equal sides and $0$ on the third side. The existence of $v$ is assured by lemma $\ref{particularbarrier}$.
We shall show that $u$ extends by continuity to $f$ along $C$. Let $P\in\partial\Omega$, fix $\epsilon>0$ and let $v$ be a barrier at $P$ defined in a triangle $V$ as described above. As $f_n$ is continuous at $P$ then $\partial\Omega$ contains a neighborhood of $P$ on which \begin{center}$f_n<f+\epsilon.$\end{center}The continuity of $f$ at $P$ allows us to assume that in this neighborhood\begin{center}$f_n<f(P)+2\epsilon$\end{center} and hence in this neighborhood \begin{center}$f_n<v+f(P)+2\epsilon.$\end{center}
We choose $K$ such that \begin{center}$\displaystyle \sup_{\partial V\cap\Omega}(u_n)<K+f(P)$\end{center}for the maximum principle would then imply the following inequality\begin{center}$u_n<v+f(P)+2\epsilon$ in $V\cap\Omega$\end{center}
Taking $n\rightarrow\infty$ implies that\begin{center}$u(x)\leq v(x)+f(P)+2\epsilon$ in $V\cap\Omega$.\end{center}
By a similar argument we obtain the inequality\begin{center}$u(x)\geq w(x)+f(P)-2\epsilon$ in $V\cap\Omega$,\end{center}where $w$ is the barrier in the triangle $V$, chosen as for $v$ above, except that $w$ takes values $-K$ on the equal sides and $0$ on the third side. The constant $K$ is chosen such that \begin{center}$\displaystyle \inf_{\partial V\cap\Omega}(u_n)>-K+f(P).$\end{center}
Taking $\epsilon\rightarrow0$ and $x\rightarrow P$, we get that $\displaystyle\lim_{x\rightarrow P}u(x)=f(P)$ and the proof is completed.\end{proof}
\section{The Conjugate Function}
Let $u$ be a solution of the minimal surface equation in a simply connected domain $\Omega$. The equation 
\begin{equation*}\displaystyle
\rm{div}_{\mathbb{R}^2}\left(\frac{\lambda\alpha}{W}\partial_x+\frac{\lambda\beta}{W}\partial_y\right)=0
\end{equation*}
amounts to the fact that the differential\begin{center}$\displaystyle \omega=\frac{-\lambda\beta}{W}dx+\frac{\lambda\alpha}{W}dy$\end{center}is exact in $\Omega$. We may then consider the function $\psi$ defined in $\Omega$, such that $d\psi=\omega$, and we shall call it the conjugate function of $u$. The gradient of $\psi$, for the $\mathbb{H}^2$-metric, is\begin{center}
$\displaystyle D\psi=\frac{-\beta}{\lambda W}\partial_x+\frac{\alpha}{\lambda W}\partial_y$
\end{center} 
and \begin{center}
$\displaystyle |D\psi|_{\mathbb{H}^2}=\frac{\sqrt{\alpha^2+\beta^2}}{W}<1,$
\end{center}
it follows that $\psi$ is Lipschitz continuous and hence extends continuously to the closure of $\Omega$ and hence $d\psi$ may be integrated along boundary arcs of $\Omega$ regardless of the boundary values of $u$. The following is obvious
\begin{lem}\label{triviallemmaondpsi}
Let $u$ be a solution of the minimal surface equation in a bounded domain $\Omega\subset\mathbb{H}^2$ and $C$ a piecewise smooth curve lying in the closure of $\Omega$. Then, \begin{center}
$\displaystyle|\int_Cd\psi|\leqslant|C|,$
\end{center} 
where $|C|$ denotes the $\mathbb{H}^2$-length of $C$.\\Moreover, if $C$ is a simple closed curve then\begin{center}
$\displaystyle\int_Cd\psi=0.$
\end{center}
\end{lem}
We remark that if $C$ lies in $\Omega$, the fact that $ |D\psi|<1$ implies that\begin{center}
$\displaystyle|\int_Cd\psi|<|C|.$
\end{center}
We show next that this will be the case when $C$ is a convex arc of the boundary of $\Omega$, provided that $u$ is continuous there.
\begin{lem}\label{strictintineq}
Let $u$ be a solution of the minimal surface equation in a domain $\Omega$ and $C$ a convex arc of the boundary of $\Omega$. If $u$ is continuous on $C$ then\begin{center}
$\displaystyle|\int_{C}d\psi|<|C|.$
\end{center}\end{lem}
\begin{proof} It is clearly enough to prove the result for a sub-arc of $C$; this allows us to assume, without loss of generality, that $\Omega$ is convex with $u$ continuous on its boundary. Let $C'$ denote the open sub-arc of the boundary which is complementary to $C$ and let $a$ be a real constant. The minimal surface equation admits a solution $u^*$ which is equal to $u$ on $C'$ and $u+a$ on $C$, as guaranteed by the above results.\\We set\begin{center}$\displaystyle\tilde{u}=u^*-u$ and $\widetilde{\psi}=\psi^*-\psi.$\end{center}
Observe that $\tilde{u}_x=-\lambda(\alpha^*-\alpha)$ and $\tilde{u}_y=-\lambda(\beta^*-\beta)$, then integration by parts and a standard "approximation" at the end-points of $C$ show that \begin{align*}
\displaystyle \int_{\partial\Omega}\tilde{u}d\widetilde{\psi}&=-\int_\Omega\Big\lbrack \tilde{u}_x\Big(\frac{\lambda\alpha}{W}-\frac{\lambda\alpha^*}{W^*}\Big)+\tilde{u}_y\Big(\frac{\lambda\beta}{W}-\frac{\lambda\beta^*}{W^*}\Big)\Big\rbrack dxdy\\
&=-\int_{\Omega}\lambda^2(\beta-\beta^*)\Big(\frac{\beta}{W}-\frac{\beta^*}{W^*}\Big)+\lambda^2(\alpha-\alpha^*)\Big(\frac{\alpha}{W}-\frac{\alpha^*}{W^*}\Big)dxdy\\
&= -\int_{\Omega}\Big\langle W\eta-W^*\eta^*,\eta-\eta^*\Big\rangle_{\psl}dA_{\mathbb{H}^2}\\
&= -\int_{\Omega}\frac{(W+W^*)}{2}(\eta-\eta^*)^2dA_{\mathbb{H}^2},
\end{align*}where $\alpha^*$, $\beta^*$, $W^*$ and $\eta^*$ are defined in terms the partial derivatives of $u^*$ in the same fashion we defined $\alpha$, $\beta$, $W$ and $\eta$ in terms of the partial derivatives of $u$. The field $\eta^*$ is normal to the graph of $u^*$.\\
The above computation then implies that \begin{center} $\displaystyle a\int_C{d\tilde{\psi}}<0.$\end{center}
Using the fact that\begin{center}$\displaystyle|\int_Cd\psi^*|\leqslant|C|$\end{center}and giving $a$ the values $\pm1$ complete the proof.\end{proof}
\begin{lem}\label{slnbehongeo}
Let $\Omega$ be a domain in $\mathbb{H}^2$ whose boundary contains a geodesic segment $\Gamma$. Suppose that $\partial\Omega$ is oriented so that the orientation on $\Gamma$ coincides with that induced by the outward pointing normal to $\Gamma$. If $u$ is a solution of the minimal surface equation  in $\Omega$ assuming boundary value plus infinity on $\Gamma$ then\begin{center}
$\displaystyle\int_\Gamma d\psi=|\Gamma|.$
\end{center} \end{lem}
\begin{proof} We consider the half plane model for the hyperbolic plane. We can suppose that $\Omega\subset\{x<0,y>0\}$ and that $\Gamma$ is a segment of the geodesic $\{x=0,y>0\}$ of $\mathbb{H}^2$.\\
We remark that the $\mathbb{H}^2$-gradient of $\psi$, the conjugate function of $u$, is \begin{center}$D\psi=Rot_{\frac{\pi}{2}}d\pi(\eta)$,\end{center} where $\eta$ is the upwards pointing unit normal to the graph $\Sigma$ of $u$, and we show that $\eta$ extends continuously to the boundary segment $\Gamma$.\\
We think of $\psl$ as a subset of $\mathbb{R}^3$ and we choose a sequence $(p_n)$ of points with constant ordinates in $\Omega$ which converges to an interior point $p$ of $\Gamma$. We set $\mu_n=d(p,p_n)$ and $q_n=(p_n,u(p_n))$ and we consider the affine transformations $\displaystyle h_n(X)=\frac{1}{\sqrt{\mu_n}}(X-q_n)$ on $\mathbb{R}^3$.\\
Let $\Sigma_n=h_n(\Sigma)$ and note that $0\in\Sigma_n$,  for all $n$, and that the normal $\eta_n$ to $\Sigma_n$ at the origin is the same as that of $\Sigma$ at the point $q_n$. It is then enough to show $\eta_n(0)$ admits a limit as $n\to\infty$ and define $\eta(p)$ as this limit.\\
We admit for now that the sequence $(A_n)$, $A_n$ the second fundamental form of $\Sigma_n$ for the euclidean metric, is uniformly bounded in a neighborhood of the origin, a claim we will prove below. Hence the sequence $\Sigma_n$ converges on this neighborhood, up to a subsequence. As $\Sigma_n$ is contained in $\{x\leqslant \sqrt{\mu_n}\}$ and asymptotic to the plane $\{x=\sqrt{\mu_n}\}$, the limit surface will be tangent to the plane $\{x=0\}$ at the origin.\vspace{.5cm}\\The sequence $(N_n)$, $N_n$ the normal to $\Sigma_n$ for the euclidean metric at the origin, therefore converges to $\partial_x$. However the equality \begin{center}$\displaystyle \eta_n=\frac{G^{-1}N_n}{\sqrt{\langle G^{-1}N_n,N_n\rangle}_{\mathbb{R}^3}},$\end{center}where the matrix $G$ is such that $\langle X,Y\rangle_{\psl}=\langle GX,Y\rangle_{\mathbb{R}^3}$, implies that $\eta_n$ is also convergent and this proves our claim that $\eta$ extends by continuity to the interior of $\Gamma$.\\
The facts that at interior points of $\Gamma$\begin{align*}\displaystyle \langle \eta,E_3\rangle_{\psl}&=\lim\langle\eta_n,E_3\rangle_{\psl}\\&=\frac{\langle \partial_x,E_3\rangle_{\mathbb{R}^3}}{\sqrt{\langle G^{-1}\partial_x,\partial_x\rangle_{\mathbb{R}^3}}}\\&=0\end{align*}and\begin{align*}\displaystyle \langle d\pi(\eta),e_2\rangle_{\mathbb{H}^2}&=\langle\eta,E_2\rangle_{\psl}\\&=\frac{\langle \partial_x,E_2\rangle_{\mathbb{R}^3}}{\sqrt{\langle G^{-1}\partial_x,\partial_x\rangle_{\mathbb{R}^3}}}\\&=0\end{align*}imply that the extension of $\eta$ to the boundary is such that\begin{center}$\langle d\pi(\eta),e_1\rangle_{\mathbb{H}^2}=-1,$\end{center}  $e_1$ being also the outwards pointing normal to $\Gamma$.\\
Now as $Rot_{\frac{\pi}{2}}$ preserves the metric on tangent spaces of $\mathbb{H}^2$ and as $\Gamma$ is oriented by $e_1$ we obtain, \begin{align*}\displaystyle\int_{\Gamma}d\psi&=-\int_{\Gamma}\langle D\psi,e_2\rangle_{\mathbb{H}^2}ds\\&=-\int_{\Gamma}\langle Rot_{\frac{\pi}{2}}d\pi(\eta),e_2\rangle_{\mathbb{H}^2}ds\\&=-\int_{\Gamma}\langle d\pi(\eta),e_1\rangle_{\mathbb{H}^2}ds\\&=|\Gamma|.\end{align*}
To complete the proof we now estimate the second fundamental form $A_n$ of $\Sigma_n$. Since $u\to\infty$ when $p\to\Gamma$ we may choose discs $\mathcal{D}(q_n,R)$ centered at $q_n$ in $\Sigma$ with intrinsic radius $R$ independent of n, and since minimal graphs in $\psl$ are stable (see Remark $\ref{stablegraphs}$), Schoen's curvature estimate implies that \begin{center}$|\widetilde{A}|\leqslant C$ in $\mathcal{D}(q_n,\frac{R}{2}),$\end{center}where $\widetilde{A}$ is the second fundamental form of $\Sigma$ for the $\psl$ metric and $C$ is an absolute constant.\\
However, if $N$ and $A$ denote the normal and the second fundamental form of $\Sigma$ with respect to the euclidean metric we have\begin{align*}\widetilde{A}(X,Y)&=\langle \overline{\nabla}_XY,\eta\rangle_{\psl}\\&=\frac{\langle \overline{\nabla}_XY,N\rangle_{\mathbb{R}^3}}{\sqrt{\langle G^{-1}N,N\rangle_{\mathbb{R}^3}}}\\&=\frac{1}{\sqrt{\langle G^{-1}N,N\rangle_{\mathbb{R}^3}}}\Big(\langle \nabla_XY,N\rangle_{\mathbb{R}^3}+\langle \overline{\nabla}_XY-\nabla_XY,N\rangle_{\mathbb{R}^3}\Big),\end{align*} 
where $\nabla$ is the Levi-Cevita connection of $\Sigma$ for the Euclidean metric. Then $\widetilde{A}$ controls $A$ as follows\begin{align*}A(X,Y)&\leqslant \sqrt{\langle G^{-1}N,N\rangle_{\mathbb{R}^3}}\widetilde{A}(X,Y)-\langle \overline{\nabla}_XY-\nabla_XY,N\rangle_{\mathbb{R}^3}.\end{align*} The tensor $ \overline{\nabla}_XY-\nabla_XY$ can be easily seen to be controlled by $\|X\|$, $\|Y\|$ and the Christofel symbols of $\psl$ which shows that $|A|$ is bounded in a neighborhood of $q_n$. Then $\widetilde{A_n}$, the second fundamental form of $\Sigma_n$ with respect to $\psl$ metric, is bounded by $\displaystyle C\sqrt{\mu_n}$ in the disc $\mathcal{D}(0,\frac{R}{2\sqrt{\mu_n}})$. In a similar fashion, one obtains the following estimates \begin{align*}A_n(X,Y)&\leqslant \sqrt{\langle G^{-1}N_n,N_n\rangle}\widetilde{A}_n(X,Y)-\langle \overline{\nabla}_XY-\nabla_XY,N_n\rangle.\end{align*}which imply that $(A_n)_n$ is uniformly bounded in a neighborhood of the origin and the proof is completed.\end{proof}
\begin{lem}\label{behaviorseqconjug}Let $\Omega$ be a domain in $\mathbb{H}^2$ as in Lemma $\ref{slnbehongeo}$ and let $(u_n)$ be a sequence of solutions of $(1)$ in $\Omega$. Assume that each $(u_n)$ is continuous in $\Omega\cup\Gamma$ and that $(u_n)$ diverges uniformly to infinity on compact subsets of $\Gamma$ while remaining uniformly bounded on compact subsets of $\Omega$. Then\begin{center}$\displaystyle\lim_{n\rightarrow\infty}\int_{\Gamma}d\psi_n=\left|\Gamma\right|.$\end{center}On the other hand, if the sequence diverges uniformly to infinity on compact subsets of $\Omega$ while remaining uniformly bounded on compact subsets of $\Gamma$, then\begin{center}$\displaystyle\lim_{n\rightarrow\infty}\int_{\Gamma}d\psi_n=-\left|\Gamma\right|.$\end{center}
\end{lem}
\begin{proof} We follow the same lines of proof as in Lemma $\ref{slnbehongeo}$ except that we choose the points $q_n=(p_n,u_n(p_n))$ instead, where $(p_n)$ is in $\Omega$ and converges to an interior point of $\Gamma$. We consider the surfaces $\mathcal{S}_n=h_n(\Sigma_n)$ in $\mathbb{R}^3$, where $\Sigma_n$ is the graph of $u_n$ and $h_n$ is as defined in the proof of Lemma $\ref{slnbehongeo}$, for purposes similar to those in that proof. We let $A_n$ denote the second fundamental form of $\mathcal{S}_n$ and $\eta_n$ the normal to $\mathcal{S}_n$ at the origin, which is the same as that to $\Sigma_n$ at $q_n$.\\
The facts that $u_n$ is continuous in $\Omega\cup\Gamma$ and that the sequence $(u_n)$ diverges uniformly on compacts of $\Gamma$, allow us to choose discs centered at $q_n$ on $\Sigma_n$, of radius $R$ independent of $n$, as in the proof of Lemma $\ref{slnbehongeo}$.\\
Moreover, we note that the sequence $(u_n)$ converges to a solution in $\Omega$ with $u$ taking the value $+\infty$ on $\Gamma$. This fact together with Schoen's curvature estimate for each $\Sigma_n$, imply in a similar way as in the proof of Lemma $\ref{slnbehongeo}$ that $(A_n)_n$ is uniformly bounded in a $D(0,R)$. The sequence $(\eta_n)$ can then be proved to converge to a horizontal vector $\eta$ along $\Gamma$ with \begin{center}$\left\langle d\pi(\eta),e_1\right\rangle=-1$\end{center} and then \begin{center}$\displaystyle\lim_{n\rightarrow\infty}\int_{\Gamma}d\psi_n=-\int_{\Gamma}{\left\langle d\pi(\eta),e_1\right\rangle} ds=\left|\Gamma\right|$\end{center}To prove the second part of the lemma we make the obvious adjustments to the proof and further details are left to the reader.\end{proof}
\section{The Monotone Convergence Theorem}
Later existence results depend on the limit behavior of monotone sequences of solutions of the minimal surface equation. In this section we develop the necessary tools to deal with these sequences. These are similar, as well as the last two sections above, to the results in \cite{js}.
\begin{lem}[Straight Line Lemma]\label{stlinelem}
Let $\Omega\subset\mathbb{H}^2$ be a bounded domain whose boundary consists of a geodesic segment $\gamma$ and an arc $C$, with $\Omega$ lying on one side of $\gamma$. Then for any compact $K\subset\Omega$ there exists $N$, depending only on the distance from $K$ to $\gamma$, such that\begin{center}
$m-N\leqslant u\leqslant M+N$ in K,
\end{center}for any solution $u$ of the minimal surface equation $(\ref{mse})$ which is bounded in $\overline{\Omega}$, with $m\leqslant u\leqslant M$ on $C$.\end{lem}
\begin{proof} Let $f_1$ and $f_2$ be two isometries of $\mathbb{H}^2$ sending the positive y-axis to the geodesic $\Gamma$ which contains $\gamma$ such that the image of the quadrant $Q_1=\{(x,y)|x>0,y>0\}$ by $f_1$ will contain $\Omega$, and the image of the quadrant $Q_2=\{(x,y)|x<0,y>0\}$ will contain $\Omega$. Let $O=f_1(Q_1)=f_2(Q_2)$ and note that the minimal graphs of example $\ref{thebarrier}$ can be used to obtain a positive solution and a negative solution of the minimal surface equation in $O$, which take respectively the value +$\infty$ and $-\infty$ on $\Gamma$.  Simply, let $\tilde{f_1}$ and $\tilde{f_2}$ denote the respective lifts of $f_1$ and $f_2$ to $\psl$ and consider the images by $\tilde{f_1}$ and $\tilde{f_2}$ of the graphs in example $\ref{thebarrier}$, which correspond to $C=1$ and defined over $Q_1$ and $Q_2$ respectively. We obtain two graphs on $O$ which, up to vertical translations, have the desired properties. Assume these graphs to be those of solutions $v_1\geq0$ and $v_2\leq0$ of $(\ref{mse})$.\\Let $u$ a solution of the minimal surface equation in $\Omega$ with $m\leqslant u\leqslant M$ on $C$. Then on the boundary of $\Omega$ we shall have\begin{center}
$m+v_2\leqslant u\leqslant M+v_1$.
\end{center}The maximum principle then implies that the inequalities hold in $\Omega$. Now, for any compact $K$ of $\Omega\cup C$, let $\displaystyle N=\max\{\max_{K}v_1,\left|\min_{K}v_2\right|\}$ which depends only on the distance from $K$ to $\gamma$. We clearly have that\begin{center}
$m-N\leqslant u\leqslant M+N$ in $K$,
\end{center}which concludes our proof.\end{proof}
\begin{rem}{\label{bndryvalues}One direct consequence of this lemma is that no solution of the minimal surface equation can take infinite values on a non-geodesic boundary arc of a convex domain. Assume to the contrary that there exists a solution $u$ of $(\ref{mse})$ in a convex domain $\Omega$ taking the value +$\infty$ (-$\infty$) on a non-geodesic open boundary arc $C$. By restricting ourselves to proper parts of $C$ we may assume $U$, the convex hull of $C$ in $\Omega$, bounded by $C$ and its end points and an open geodesic segment $\gamma$ contained in $\Omega$. We shall obtain a contradiction by showing that $u$ must be equal to +$\infty$ (-$\infty$) in $U$.\\Let $\displaystyle a=\inf_\gamma u$ ($\displaystyle =\sup_\gamma u$) which may be assumed a positive (negative) real number (if we restrict ourselves to proper parts of $C$). For each $n$, let $u_n$ be the solution of the minimal surface equation in $U$ taking the values $n$ $(-n)$ on $C$ and $a$ on $\gamma$. By the maximum principle, we have then $u_n\leqslant u$ ($u_n\geqslant u$) in $U$. Hence by the Straight Line Lemma we have that on each compact in $U\cup C$ and for each $n$, $n-N\leqslant u_n\leqslant u$ ($-n+N\geqslant u_n\leqslant u$) with $N$ independent of $n$. Letting $n\to\infty$ implies that $u$ has infinite values in $U$ which is absurd.}\end{rem}
The following two theorems are essential for studying convergence of monotone sequences of solutions.
\begin{thm}[Monotone Convergence Theorem]\label{monotonecv}Let $(u_n)$ be a monotonically increasing sequence of solutions of the minimal surface equation in a domain $\Omega$. If the sequence is bounded at a point $p\in\Omega$, then there exists a non-empty open set $U\subset\Omega$ such that the sequence $(u_n)$ converges to a solution in $U$, and diverges to infinity on the complement of $U$. The convergence is uniform on compacts of $U$, and the divergence is uniform on compacts of $V=\Omega-U$.\end{thm}\rm{}
\begin{proof}  Assume that $\left|u_0\right|\leq c$ near $p$, and consider the sequence of non-negative solutions $(v_n)$ such that $v_n=u_n+c$. Hence, each $W_n(p)\leq C_n$, where $C_n$ is the constant given by theorem $\ref{gradientestimates}$. Then theorem $\ref{uniformbounds}$ implies that each $W_n$ is bounded in a disc centered at $p$ and of radius $R_n$, with $R_n$ depending on $u_n(p)$, $d(p,\partial\Omega)$ and bounds of $\lambda$ and its derivatives. As $(u_n(p))$ is bounded, then we can find a disc $D$ centered at $p$ on which $(W_n)$ is uniformly bounded. The mean value theorem then implies that $(u_n)$ is then uniformly bounded in this disc. The compactness principle therefore implies that $(u_n)$ has a convergent subsequence and as $(u_n)$ is monotone it converges on this disc. The compactness principle implies also that the limit is a solution of the minimal surface equation and so $U$ is a non-empty open set.  The divergence is uniform on compacts of $V$ as the sequence is monotonically increasing.\end{proof}
The divergence set $V$ is by no means arbitrary, it has a very particular geometric structure. We resume the properties of $V$ in the following
\begin{thm}[Divergence Set Structure Theorem]
\label{divergencesetstructure}Let $(u_n)$ be a monotonically increasing sequence of solutions in $\Omega$. If the divergence set $V\neq\emptyset$, then int$(V)\neq\emptyset$, and $\partial V$ is composed of non-intersecting geodesic segments of $\Omega$ and possibly parts of $\partial\Omega$. Moreover, no two interior geodesic segments of $\partial V$ can have a common end point at a convex corner of $V$, nor any component of $V$ consist only of a geodesic segment of $\Omega$.\\
Furthermore, if $\Omega$ is bounded in part  by a convex arc $C$ with each $u_n$ continuous in $\Omega\cup C$ and $(u_n)$ either diverges to infinity on $C$ or remains uniformly bounded on compacts of $C$, then no interior geodesic segment $\Gamma$ forming part of the boundary of $V$ can terminate at an interior point of $C$.\end{thm}
For the proof of this theorem, one can employ the lemmas of section $6$ above in ways similar to those in the proofs of Lemma $5$ and Lemma $6$ in \cite{js}. In fact, Remark $\ref{bndryvalues}$ above implies that, if $V\neq\emptyset$, $\partial V$ consists of non intersecting geodesic segments of $\Omega$ and possibly parts of the boundary of $\Omega$. To prove that no component of $\partial V$ is only a geodesic segment $T$ of $\Omega$, one applies Lemma $\ref{behaviorseqconjug}$ above to $\Omega_1$ and $\Omega_2$, the components of $\Omega$ on either side of $T$. A contradiction is obtained since in $\Omega_1$, say, one obtains \begin{center}$\displaystyle\lim_{n\rightarrow\infty}\int_{T}{d\psi_n}=\left|T\right|$\end{center}
and in $\Omega_2$ one obtains \begin{center}$\displaystyle\lim_{n\rightarrow\infty}\int_{T}{d\psi_n}=-\left|T\right|.$\end{center}
To see that no interior geodesic segments of $\partial V$ can have a common end point, we notice that Remark $\ref{bndryvalues}$ above implies that such a point must be in $\partial\Omega$. We suppose then $\partial V$ admits two geodesic segments $T_1$ and $T_2$ in $\Omega$ with a common end point $Q$ in $\partial\Omega$ and we choose two points $Q_1$ and $Q_2$ on $T_1$ and $T_2$ respectively, so that the open geodesic triangle $\Delta$, with vertices $Q$, $Q_1$ and $Q_2$, lie in $\Omega$. By Lemma $\ref{triviallemmaondpsi}$ above,
\begin{center}$\displaystyle\int_{QQ_1}{d\psi_n}+\int_{Q_1Q_2}{d\psi_n}+\int_{Q_2Q}{d\psi_n}=0$.\end{center} The triangle may lie in $U$ or $V$, since no component of $\partial V$ is only a geodesic segment. In the former case, one applies Lemma $\ref{behaviorseqconjug}$ above to obtain
\begin{center}$\displaystyle\lim_{n\rightarrow\infty}\int_{QQ_1}{d\psi_n}=\left|QQ_1\right|$ and $\displaystyle\lim_{n\rightarrow\infty}\int_{QQ_2}{d\psi_n}=\left|QQ_2\right|$\end{center}assuming that $QQ_1Q_2$ determines the positive orientation of $\Delta$. However, Lemma $\ref{triviallemmaondpsi}$ implies \begin{center}$\displaystyle\left|\int_{Q_1Q_2}d\psi_n\right|\leq\left|Q_1Q_2\right|$\end{center}which is a contradiction with the triangle inequality in $\mathbb{H}^2$. If $\Delta$ lies in $V$ one obtains a similar contradiction by applying the second part of Lemma $\ref{behaviorseqconjug}$.
To prove the second part of the theorem , we notice that if $C$ is not geodesic, Lemma $\ref{stlinelem}$ implies that on compacts in the convex hull of $C$\begin{center}$\displaystyle\min_{C}(u_n)-N\leq u_n\leq\max_{C}(u_n)+N$\end{center}with $N$ independent of $n$, and the proof of the claim is immediate since the above inequality implies that the interior of the convex hull of $C$ lies either in $U$ or in $V$.\\ We then assume that $C$ is geodesic, and that $\Gamma$ terminates at an interior point $Q$ of $C$. Suppose first that the sequence diverges on $C$. Let $P$ be a point of $\Gamma$, and choose a point $R$ on $C$ such that the geodesic segment $RP$ lies in $U$. The results we have proved in the first part of the theorem allow this choice. We apply Lemma $\ref{triviallemmaondpsi}$ and Lemma $\ref{behaviorseqconjug}$ in the triangle $QPR$ in a fashion similar to that in the triangle $\Delta$, to obtain a contradiction with the triangle inequality in $\mathbb{H}^2$. In case the sequence remains uniformly bounded on compacts of $C$, a similar contradiction results by choosing the segment $RP$ in $V$.$\hfill\Box$
\section{A Jenkins-Serrin Type Theorem}
This is Theorem $\ref{jsmainthm}$ stated in the introduction. We note that a section  $s$ of $\pi:\psl\to\mathbb{H}^2$ takes the value $+\infty$ ($-\infty$ resp.) on a geodesic segment $A_i$ ($B_j$ resp.) if the image by $s$ of each geodesic $t\to\gamma(t)$ of $\Omega$ ending at $A_i$ ($B_j$ resp.) gets out of every compact and if $\langle \gamma'(t),\xi\rangle>0 (<0, resp.)$, where $\xi=\partial_z$ in our model of $\psl$.

As was remarked above, having fixed the model for $\psl$ the existence of the section $s$ on $\Omega$ with the prescribed boundary data is equivalent to the existence of a real function $u$ defined in $\Omega$ with corresponding data on the boundary. The function $u$ is constructed as a limit of monotone sequence of solutions of the minimal surface equation whose behavior is studied using the monotone convergence theorem, the divergence set structure theorem and the properties of the differential $d\psi$ corresponding to $u$. Once the convergence is established, we need to show that the limit will assume the appropriate boundary values. This will be assured by the Boundary Values Lemma below.

We proceed to prove the existence of particular solutions of $(\ref{mse})$ which will be used as barriers in the proof of the Boundary Values Lemma.
\begin{lem}\label{bdryvaluesbarrier} Let $\mathcal{P}$ be a convex quadrilateral in $\mathbb{H}^2$, formed by geodesic segments $A_1$, $A_2$, $C_1$ and $C_2$ such that $A_1\cap A_2=\emptyset$ and $|A_1|+|A_2|<|C_1|+|C_2|$. Then there exists a solution of $(\ref{mse})$ in $\mathcal{P}$ which takes the boundary values $+\infty$ on $A_1\cup\ A_2$ and non-negative values on $C_1\cup C_2$.\end{lem}
\begin{proof} Let $\mathcal{P}$ be a convex quadrilateral in $\mathbb{H}^2$, formed by geodesic segments $A_1$, $A_2$, $C_1$ and $C_2$ such that $A_1\cap A_2=\emptyset$ and $|A_1|+|A_2|<|C_1|+|C_2|$. Let $u_n$ be the solution of the minimal surface equation in $\mathcal{P}$ taking boundary values $n$ on each $A_i$ and 0 on each $C_i$. The sequence $u_n$ is seen to converge to a solution $u$ in $\mathcal{P}$ as follows. Let $\mathcal{V}$ denote the divergence set and remark that either the interior of $\mathcal{V}$ is equal to that of $\mathcal{P}$, or otherwise by Theorem $\ref{divergencesetstructure}$ an interior geodesic segment bounding $\mathcal{V}$ must have its endpoints from amongst those of the $A_i$'s.\vspace{.25cm}\\
The interior of $\mathcal{V}$ cannot be equal to that of $\mathcal{P}$ for otherwise:\begin{center}
$\displaystyle\int_{A_1\cup\ A_2}{d\psi_n}+\int_{C_1\cup\ C_2}{d\psi_n}=0$ \end{center}
and then one takes the limit as $n\rightarrow+\infty$ and uses Lemma $\ref{behaviorseqconjug}$ and Lemma $\ref{triviallemmaondpsi}$ to obtain \begin{center}
$\displaystyle\int_{C_1\cup\ C_2}{d\psi_n}=-\big(|C_1|+|C_2|\big)$ and $\displaystyle\int_{A_1\cup\ A_2}{d\psi_n}\leq |A_1|+|A_2|$.\end{center}
This implies that $|A_1|+|A_2|\geq|C_1|+|C_2|$ which is not true.\vspace{.25cm}\\
Thus assume that $\mathcal{V}$ is non-empty and bounded by a geodesic triangle $\Delta$ whose vertices are endpoints of the $A_i$'s. Let $\delta$ denote the perimeter of $\Delta$. One would obtain \begin{center}
$\displaystyle\int_{\Delta-A_i}{d\psi_n}+\int_{A_i}{d\psi_n}=0$.\end{center}
Again passing to the limit and using Lemma $\ref{behaviorseqconjug}$ and Lemma $\ref{triviallemmaondpsi}$ the following holds\begin{center}
$\displaystyle\int_{\Delta-A_i}{d\psi_n}=-(\delta-|A_i|)$ and $\displaystyle \int_{A_i}{d\psi_n}\leq |A_i|$,\end{center}
which leads to a contradiction with the triangle inequality.\vspace{.25cm}\\
Therefore, $\mathcal{V}=\emptyset$ and the sequence $(u_n)$ converges on compact sets of $\mathcal{P}$ to a solution of ($\ref{mse}$). We note that since $(u_n)$ is increasing (by the maximum principle), $u$ takes the value $+\infty$ on the segments $A_i$. Although at this point we do not know yet that $u=0$ on the
$C_i$'s, a fact which we will be able to prove later, we remark that $u\geq0$ on each $C_i$.\end{proof}
\begin{lem}\label{baronrec} Let $\mathcal{P}$ be a convex quadrilateral in $\mathbb{H}^2$ formed by geodesic segments $A_1$, $A_2$, $C_1$ and $C_2$ such that $A_1\cap A_2=\emptyset$. If $|A_1|+|A_2|<|C_1|+|C_2|$, then there exists a solution $v$ of the minimal surface equation in $\mathcal{P}$ taking the boundary values $+\infty$ on $A_1\cup\ A_2$ and has bounded values on $C_1\cup C_2$.\end{lem}
\begin{proof} Let $\widetilde{C_i}$ be a horizontal lift of $C_i$ to $\psl$ and let $u_n$ be the solution of the minimal surface equation in $\mathcal{P}$ taking boundary values $n$ on each $A_i$ and boundary values given by $\widetilde{C_i}$ on $C_i$. One may translate vertically each of the $\widetilde{C_i}$, so that each $u_n\geq0$. The sequence $u_n$ is increasing and converges to a solution $u$ in $\mathcal{P}$ by arguments similar to those in Lemma $\ref{bdryvaluesbarrier}$.\\
The limit $u$ takes the boundary value $+\infty$ on each $A_i$ as the sequence $(u_n)$ is increasing. The boundary values of $u$ on $C_i$ are given by $\widetilde{C_i}$ and this follows by standard barrier techniques, as in the proof of Proposition $\ref{radogen}$, once we show the sequence $(u_n)$ to be uniformly bounded near each point of $C_i$. We complete the proof by showing this last point.\\As $\widetilde{C_i}$ is a horizontal geodesic, observations in \cite{Uwe} ensures that the graph of $u_n$ extends by symmetry about $\widetilde{C_i}$ to a graph (a graph since otherwise by the maximum principle, the surface obtained by symmetry would coincide with the cylinder $\pi^{-1}(C_i)$).\\
Let $p$ be a point of $C_i$ and choose a sufficiently small geodesic rectangle $R$ as in Lemma $\ref{bdryvaluesbarrier}$, which has two of its sides orthogonal to $C_i$ and which contains $p$ in its interior. Let $v$ denote the solution of ($\ref{mse}$) in $R$ taking the values $+\infty$ on the sides orthogonal to $C_i$ and non-negative values on the other two sides, say $S_j$, $1\leq j\leq2$. The existence of $v$ is assured by Lemma $\ref{bdryvaluesbarrier}$. The maximum principle then implies that for each $n$, $u_n\leq v+M$ over $R$, where $\displaystyle M=\sup{\left|u_n\right|}$ and the supremum is taken over the $S_j$'s. One considers a small neighborhood of $p$ in $R$, and the preceding inequality proves that $u_n$ is bounded around $p$.\end{proof}

\begin{lem}[Boundary Values Lemma]Let $\Omega$ be a domain and $C$ a compact convex arc in its boundary. Let $(u_n)$ be a sequence of solutions of the minimal surface equation, which converges uniformly on compacts of $\Omega$ to a solution $u$. Suppose that, on the one hand, each $(u_n)$ is continuous in $\Omega\cup C$ and that the boundary values converge uniformly on compacts of $C$ to a limit function $f$. Then $u$ is continuous in $\Omega\cup C$ and takes the values $f$ on $C$. If $C$, on the other hand, were a geodesic segment where the boundary values diverge uniformly to infinity, then $u$ will take on the boundary value infinity on $C$.\end{lem}
\begin{proof} For the first part where the boundary values of $(u_n)$ converge uniformly on compact subsets of $C$ it suffices to show the sequence $(u_n)$ uniformly bounded in the neighborhood of any interior point of $C$ and then employ a standard argument of the theory of barriers (again similar to that in the proof of proposition $\ref{radogen}$ above). If $C$ is not geodesic then the result follows by the Straight Line Lemma. In case $C$ is a geodesic segment, the preceding lemma furnishes the ingredient necessary to show the required boundedness of $(u_n)$ near interior points of $C$ in a fashion similar to that of Lemma 7 in \cite{js} or the corresponding Boundary Values Lemma in \cite{rc}.\\The part where $C$ is geodesic and $(u_n)$ taking infinite values there can be proved in a fashion similar to that of Lemma $8$ in \cite{js}. However, to prove $(u_n)$ bounded from below as is done in \cite{js} we may prove a lemma similar to Lemma $\ref{baronrec}$ above except that the solution takes values $-\infty$ on the sides $A_i$. Then we follow the same lines of proof of the Boundary Values Lemma in \cite{rc}.\end{proof}
\begin{rem}Let $\mathcal{P}$ be a geodesic rectangle as in Lemma $\ref{baronrec}$. In order to prove the existence of a solution of ($\ref{mse}$) in $\mathcal{P}$ taking bounded values on the $C_i$'s and values $-\infty$ on the $A_i$'s, one can proceed as follows.
Let $r$ denote the reflection of $\mathbb{H}^2$ in $A_1$, and $\tilde{r}$ its lift to $\psl$. Consider the image $\mathcal{P}'$ of $\mathcal{P}$ by $r$, and find by Lemma $\ref{baronrec}$ a solution $v$ in $\mathcal{P}'$ taking the values $+\infty$ on $A_1$ and $r(A_2)$, and bounded values on each $r(C_i)$. The image by $\tilde{r}$ of the graph of $v$ is the graph of a solution $u$ of ($\ref{mse}$) over $\mathcal{P}$, which takes the sought boundary values.
\end{rem}
Having developed the necessary machinery in this paper, the existence part of Theorem $\ref{jsmainthm}$ could be proved following the same lines of proof in \cite{js} and \cite{nelli}. To see that the conditions in Theorem $\ref{jsmainthm}$ are necessary, we let $u$ be a solution of the minimal surface equation in a domain $\Omega$ and we consider a polygon $\mathcal{P}$, such that $\Omega$ and $\mathcal{P}$ are as described in that theorem. By Lemma $\ref{triviallemmaondpsi}$ above \begin{center}$\displaystyle\int_{\mathcal{P}-\{A_i\in\mathcal{P}\}}{d\psi}+\sum_{A_i\in\mathcal{P}}{\int_{A_i}{d\psi}}=0$,\end{center}with $\mathcal{P}$ oriented by the outward pointing normal. Lemma $\ref{slnbehongeo}$ implies that\begin{center}$\displaystyle\sum_{A_i\in\mathcal{P}}{\int_{A_i}{d\psi}}=\alpha$,\end{center}and if  $\mathcal{P}\neq\partial\Omega$ then Lemma $\ref{strictintineq}$ implies that \begin{center}$\displaystyle\int_{\mathcal{P}-\{A_i\in\mathcal{P}\}}{d\psi}<\gamma-\alpha$.\end{center}If $\mathcal{P}=\partial\Omega$, which is possible only if the family $\{C_i\}=\emptyset$, then again by Lemma $\ref{slnbehongeo}$, one would obtain
\begin{center}$\displaystyle\int_{\mathcal{P}-\{A_i\in\mathcal{P}\}}{d\psi}=-\beta$.\end{center}
This argument shows that the conditions $2\alpha<\gamma$ for all possible polygons $\mathcal{P}\neq\partial\Omega$ chosen as in Theorem $\ref{jsmainthm}$, and that $\alpha=\beta$ when $\mathcal{P}=\partial\Omega$ are necessary. A similar argument shows that the conditions on the segments $B_i$ are necessary as well.

To show that the conditions of Theorem $\ref{jsmainthm}$ are sufficient, we employ the Monotone Convergence Theorem, the Divergence Structure Theorem and the lemmas of sections 6 through 8 in the same fashion as in \cite{js} or \cite{nelli}. We furnish only a sketch of the proof and we refer the reader to section 5 in \cite{js} for further details, where the constructions of solutions held in that paper carry word for word to our case. The proof can be broken down into proving existence of solutions of Dirichlet problems related to the one stated in Theorem $\ref{jsmainthm}$.\\
\textbf{Step 1}. We consider the Dirichlet problem in Theorem $\ref{jsmainthm}$, and we suppose that the family $\{B_i\}$ is empty. Assume also that the assigned data on the arcs $\{C_i\}$ is bounded below. Then the conditions $2\alpha<\gamma$ for each simple closed polygon $\mathcal{P}$ whose vertices are chosen from among the endpoints of the $A_i$'s are sufficient for the existence of a solution.\\
\textbf{Step 2}. We consider the Dirichlet problem in Theorem $\ref{jsmainthm}$, and we suppose that the family $\{C_i\}\neq\emptyset$. Then the conditions $2\alpha<\gamma$ and $2\beta<\gamma$ for each simple closed polygon $\mathcal{P}$ whose vertices are chosen from among the endpoints of the $A_i$'s and the $B_i$'s are sufficient for the existence of a solution.\\
\textbf{Step 3}. We consider the Dirichlet problem in Theorem $\ref{jsmainthm}$, and we suppose the family $\{C_i\}=\emptyset$. Then the conditions $\alpha=\beta$ when $\mathcal{P}=\partial\Omega$ and $2\alpha<\gamma$ and $2\beta<\gamma$ for each simple closed polygon $\mathcal{P}$ whose vertices are chosen from among the endpoints of the $A_i$'s and the $B_i$'s are sufficient for the existence of a solution.

To complete the proof of Theorem $\ref{jsmainthm}$, we next give a proof of the uniqueness, which is up to an additive constant when the family $\{C_i\}=\emptyset$, inspired by \cite{collin}.

\textbf{Proof of uniqueness}. Let $u_1$ and $u_2$ be two different solutions of the minimal surface equation  with the same boundary data (possibly infinite). If $\{C_i\}=\emptyset$ we suppose that $u_1-u_2$ is not a constant. Note that either of the subset of $\Omega$, $\{u_1>u_2\}$ or $\{u_1<u_2\}$ is non-empty . We suppose without loss of generality that $\{u_1>u_2\}\neq\emptyset$ and we choose $\epsilon$ small enough so that $\Omega_\epsilon=\{u_1-u_2>\epsilon\}$ is non-empty and that $\partial\Omega_\epsilon$ is regular.\\We consider the closed differential $d\Psi=d\psi_1-d\psi_2$, $\psi_1$ and $\psi_2$ the conjugate functions of $u_1$ and $u_2$ respectively, and we obtain a contradiction by showing that $\int_{\partial\Omega_\epsilon} d\Psi\neq0$.\\As $u_1$ and $u_2$ have the same boundary data, $\partial\Omega_\epsilon$ does not intersect $\cup C_i$, besides Lemma $\ref{slnbehongeo}$ implies that $d\Psi=0$ on $\cup A_i\bigcup\cup B_j$. Then the only part of $\partial\Omega_\epsilon$ which contributes to the integral $\int_{\partial\Omega_\epsilon} d\Psi$, denoted $\widetilde{\partial\Omega_\epsilon}$, is that contained in $\Omega_\epsilon$ defined by $u_1-u_2=\epsilon$. Then the  vector \begin{center}$v=Rot_\frac{\pi}{2}\Big(\nabla(u_1-u_2)\Big)=-\lambda(\beta_1-\beta_2)\partial_x+\lambda(\alpha_1-\alpha_2)\partial_y$\end{center} is tangent to $\widetilde{\partial\Omega_\epsilon}$ and the integral $\int_{\partial\Omega_\epsilon} d\Psi$ reduces to integrating $d\Psi.v$. However, a computation similar to that of lemma $\ref{strictintineq}$ shows that\begin{center} $\displaystyle d\Psi.v=\lambda^2\frac{(W_1+W_2)}{2}(\eta_1-\eta_2)^2$
\end{center}
which is a positive quantity ($\eta_i$ is the normal to the graph of $u_i$). This leads to a contradiction and the proof is completed.$\hfill\Box$
\bibliographystyle{plain}
\bibliography{mspsl}
\nocite{*}
\end{document}